\documentclass[10pt]{amsart}
\usepackage{amsmath}
\usepackage{amssymb}
\usepackage{amsfonts}
\usepackage{mathrsfs}
\usepackage{graphicx}
\usepackage{epsfig}

\DeclareMathAlphabet{\mathpzc}{OT1}{pzc}{m}{it}

\theoremstyle{plain}
\newtheorem{theo}{Theorem}
\newtheorem*{theo }{Theorem}
\newtheorem{prop}{Proposition}[section]

\newtheorem{coro}[prop]{Corollary}
\newtheorem{lemm}[prop]{Lemma}
\newtheorem{theoalph}{Theorem}

\theoremstyle{definition}
\newtheorem{defi}[prop]{Definition}

\theoremstyle{remark}
\newtheorem{rema}[prop]{Remark}

\newtheoremstyle{citing}
  {3pt}
  {3pt}
  {\itshape}
  {}
  {\bfseries}
  {.}
  {.5em}
  {\thmnote{#3}}

\theoremstyle{citing}
\newtheorem*{generic}{}

\newcommand{\partn}[1]{{\smallskip \noindent \textbf{#1.}}}

\newcommand\C{\mathbb{C}}
\newcommand\D{\mathbb{D}}

\newcommand\R{\mathbb{R}}

\newcommand\cB{\mathcal{B}}
\newcommand\cC{\mathcal{C}}

\newcommand\fD{\mathfrak{D}}
\newcommand\fm{\mathfrak{m}}
\newcommand\fF{\mathfrak{F}}

\newcommand\tB{\widetilde{B}}
\newcommand\tc{\widetilde{c}}
\newcommand\tD{\widetilde{D}}

\newcommand\tF{\widetilde{F}}

\newcommand\tK{\widetilde{K}}

\newcommand\tU{\widetilde{U}}
\newcommand\tV{\widetilde{V}}
\newcommand\tW{\widetilde{W}}

\newcommand\sB{\mathscr{B}}
\newcommand\sP{\mathscr{P}}

\newcommand\hD{\widehat{D}}

\newcommand\hU{\widehat{U}}
\newcommand\hV{\widehat{V}}
\newcommand\hW{\widehat{W}}

\newcommand\hfF{\widehat{\fF}}

\newcommand\uW{\underline{W}}

\newcommand\e{\varepsilon}
\newcommand\la{\lambda}
\newcommand\ov{\overline}

\renewcommand{\=}{ : = }

\newcommand{\axiomM}[1]{{({\rm I}$_{#1}$)}}

\DeclareMathOperator{\diam}{diam}
\DeclareMathOperator{\Jac}{Jac}

\DeclareMathOperator{\dist}{dist}
\DeclareMathOperator{\Area}{Area}

\DeclareMathOperator{\modulus}{mod}
\newcommand{\CC}{\overline{\C}}
\newcommand{\Crit}{\text{Crit}}

\newcommand{\HD}{\text{HD}}
\newcommand{\BD}{\overline{\text{BD}}}

\newcommand{\talpha}{\widetilde{\alpha}}
\newcommand{\tnu}{\widetilde{\nu}}

\newcommand{\Lip}{\text{Lip}}
\newcommand{\map}{f}
\newcommand{\CJ}{\Crit(\map) \cap J(\map)}

\newcommand{\sC}{\mathscr{C}}
\newcommand{\sCJ}{\sC}
\newcommand{\modul}{\mathbf{m}}
\newcommand{\bfm}{\mathbf{m}}
\newcommand{\measure}{\mathbf{\fm}}

\newcommand{\hmu}{\widehat{\mu}}

\newcommand{\HDhyp}{\HD_{\text{hyp}}}
\newcommand{\Jcon}{J_{\text{con}}}
\newcommand{\Exp}{\text{Exp}}

\newcommand{\TCE}{Topological Collet-Eckmann Condition}
\newcommand{\TCEC}{TCE condition}
\newcommand{\CEC}{Collet-Eckmann condition}
\newcommand{\ExpShrink}{Exp\-Shrink}
\newcommand{\Na}[1]{{\left\| #1 \right\|_{\alpha}}}

\newcommand{\axiomC}{{{\rm (C)}}}
\newcommand{\CLT}{Central Limit Theorem}

\DeclareMathOperator{\supp}{supp}

\begin{document}

\title[Statistical properties of Topological Collet-Eckmann maps]{Statistical properties of Topological Collet-Eckmann maps.}
\author[F. Przytycki]{Feliks Przytycki$^\dag$}
\author[J. Rivera-Letelier]{Juan Rivera-Letelier$^\ddag$}
\thanks{$\dag$ Partially supported by Polish KBN grant 2 P03A 03425.}
\thanks{$\ddag$ Partially supported by Fondecyt N~10040683, MeceSup UCN-0202, IMPAN and Polish KBN grant 2 P03A 03425.}
\address{$\dag$ Institute of Mathematics Polish Academy of Sciences, 
ul. \'Sniadeckich 8, 00956 Warszawa, Poland.}
\email{feliksp@impan.gov.pl}
\address{$\ddag$ Departamento de Matem\'aticas, Universidad Cat\'olica del Norte, Casilla~1280, Antofagasta, Chile.}
\email{rivera-letelier@ucn.cl}
\date{February 2006.}
\begin{abstract}
We study geometric and statistical properties of complex rational maps satisfying the \TCE.
We show that every such a rational map possesses a unique conformal probability measure of minimal exponent, and that this measure is non-atomic, ergodic and that its Hausdorff dimension is equal to the Hausdorff dimension of the Julia set.
Furthermore, we show that there is a unique invariant probability measure that is absolutely continuous with respect to this conformal measure, and we show that this measure is exponentially mixing (it has exponential decay of correlations) and that it satisfies the \CLT.

We also show that for a complex rational map $f$ the existence of such an invariant measure characterizes the \TCE, and that this measure is the unique equilibrium state with potential $ -\HD(J(f)) \ln |f'|$.
\end{abstract}

\maketitle



\section{Introduction}
We consider complex rational maps $f : \CC \to \CC$ of degree at least~$2$, viewed as dynamical systems acting on the Riemann sphere~$\CC$.
We provide a systematic approach to study geometric and statistical properties of rational maps.
For simplicity we restrict to rational maps satisfying the ``\TCE'' (TCE).
This condition is very natural and important, because it has several equivalent formulations~\cite{PRS} and because the set of (non-hyperbolic) rational maps satisfying this condition has positive Lebesgue measure in the space of all rational maps of a given degree~\cite{Asp}.
Our main results extend without change to multimodal maps of the interval satisfying the \TCEC.
\subsection{The \TCE{}}\label{si:TCE}
The \TCEC{} was originally formulated in topological terms, but it has many equivalent formulations.
For example, a rational map $f$ satisfies the \TCEC{} if and only if it is \textit{non uniformly hyperbolic}: The Lyapunov exponent of each invariant probability measure supported on the Julia set is larger than a positive constant, that is independent of the measure.
In this paper we will mainly use the following equivalent formulation of the \TCEC.

\begin{generic}[Exponential Shrinking of Components (\ExpShrink)]
There exist $\la_{\Exp}>1$ and $r_0 > 0$ such that for every $x \in J(f)$, every integer $n \ge 1$ and every connected component $W$ of $f^{-n}(B(x, r_0))$ we have
$$
\diam(W) \le \la_{\Exp}^{-n}.
$$
\end{generic}
See~\cite{PRS} for the original formulation of the \TCEC, and for other of its equivalent formulations.

The \TCEC{} is closely related to the ``\CEC'': a rational map $f$ satisfies the \textit{\CEC} if every non-repelling periodic point of~$f$ is attracting, and if for every critical value~$v$ in the Julia set $J(f)$ of~$f$ that is not mapped to a critical point under forward iteration, the derivative $|(f^n)'(v)|$ growths exponentially with~$n$.
The \CEC{} implies the \TCEC, but the reserve implication is not true.
In fact, in \cite[\S5]{PRo} there is an example of a rational map~$f$ satisfying the \TCEC{} and having a critical value~$v$ in $J(f)$ that is not mapped to a critical point under forward iteration, and such that $\liminf_{n \to \infty} \frac{1}{n} \ln |(f')^n(v)| = - \infty$.

The \CEC{} was introduced in~\cite{CE83}, in the context of unimodal maps.
It has been extensively studied for complex rational maps, see~\cite{Asp,DF,GSmCE,GSw,P2,P1,PRS, PRo,SmCE} and references therein.
In particular, M.~Aspenberg recently proved in~\cite{Asp} that in the space of all rational maps of a given degree, there is a set of positive Lebesgue measure of (non-hyperbolic) rational maps that satisfy the \CEC, and hence the \TCEC.

\subsection{Conformal measures.}\label{si:conformal measures}
A general Julia set has a fractal nature, as its Hausdorff dimension is larger than its topological dimension, see e.g.~\cite{Zd}.
In this case, a natural geometric measure on the Julia set are the conformal measures of minimal exponent: given $t > 0$, a non zero Borel measure $\mu$ is {\it conformal of exponent} $t$ for~$f$, if for every Borel subset $U$ of $\CC$ where~$f$ is injective, we have
$$
\mu(f(U)) = \int_U |f'|^t d \mu.
$$
Every rational map admits a conformal measure~\cite{Su} and the minimal exponent for which such a measure exists is equal to the ``hyperbolic dimension'' of the Julia set, see~\cite{DU, PrLyapunov} and also~\cite{McM, U, PUbook}.
For a uniformly hyperbolic rational map there is a unique conformal probability measure of minimal exponent.
This measure is equal, up to a constant factor, to the restriction to the Julia set of the Hausdorff measure of dimension equal to the Hausdorff dimension of the Julia set~\cite{Su}.
However, for certain rational maps the Hausdorff measure of the Julia set might be zero or infinity.
In other cases all the conformal measures are atomic.
See the survey article~\cite{U} for these and other results related to conformal measures.

Our first result is about the existence of a non-atomic conformal measure of minimal exponent.
\begin{theoalph}\label{t:conformal measures}
Every rational map satisfying the \TCEC{} admits a unique conformal probability measure of minimal exponent.
This measure is non-atomic, ergodic, its Hausdorff dimension is equal to the Hausdorff dimension of the Julia set and it is supported on the conical Julia set.
\end{theoalph}
We recall the definition of ``conical Julia set'' in~Appendix~\ref{a:conformal via inducing}.
This set is also called ``radial Julia set''.

The conformal measures of a rational map without recurrent critical points in the Julia set are well understood~\cite{U}.
The first result in the more subtle case when there is a recurrent critical point in the Julia set was proved by the first named author in~\cite{P2}, where it is shown that a rational map satisfying the \CEC{} and an additional ``Tsujii type'' condition admits a non-atomic conformal measure of minimal exponent.
This result was extended by J.~Graczyk and S.~Smirnov to rational maps satisfying the \CEC, and the weaker ``summability condition''~\cite{GSm}.
The methods employed in these articles breakdown for rational maps satisfying the \TCEC, as they use the growth of derivatives at critical values in an essential way.

\subsection{Absolutely continuous invariant measures and their statistical properties.}\label{si:invariant measures}
Having a non-atomic conformal measure of minimal exponent as a reference measure, it is natural to look for absolutely continuous invariant measures.
Note that if a rational map satisfies the \TCEC{} and its Julia set is the whole sphere, then the measure given by Theorem~\ref{t:conformal measures} is the spherical measure.

Recall that if $(X, \sB, \nu)$ is a probability space and if $f : X \to X$ is a measure preserving map, then the measure $\nu$ is said to be \textit{mixing}, if for every pair of square integrable functions $\varphi, \psi : X \to \R$ we have
\begin{equation*}
\int (\varphi \circ f^n) \cdot \psi d\mu - \int \varphi d\nu \int \psi d\nu
\to 0, \text{as $n \to \infty$}.
\end{equation*}
\begin{theoalph}\label{t:invariant measures}
Let $f$ be a rational map satisfying the \TCEC.
Then there is a unique~$f$ invariant probability measure that is absolutely continuous with respect to the unique conformal probability measure of minimal exponent of~$f$.
Moreover this measure is ergodic, mixing and its density with respect to the conformal measure of minimal exponent is almost everywhere bounded from below by a positive constant.
\end{theoalph}
There are several existence results of this type for rational maps with no recurrent critical points in their Julia set~\cite{U}.
For other existence results in the case of recurrent critical points, see~\cite{Asp,Bernard,GSm,P1, P2,Rees}.

We next study the statistical properties of the invariant measure given by Theorem~\ref{t:invariant measures}.
Given a measurable space $(X, \sB, \nu)$ and a measure preserving map $f : X \to X$, we will say that the measure $\nu$ is {\it exponentially mixing} or that it has {\it exponential decay of correlations}, if there are constants $C > 0$ and $\rho \in (0, 1)$ such that every bounded measurable function $\varphi$ and every Lipschitz continuous function $\psi$, we have
\begin{equation}\label{e:correlation bound}
\left| \int (\varphi \circ f^n) \cdot \psi d\nu - \int \varphi d\nu \int \psi d\nu \right|
\le C \left( \sup_{z \in J(f)} |\varphi(z)| \right) \| \psi \|_{\Lip} \cdot \rho^n.
\end{equation}
Here $\| \psi \|_{\Lip} = \sup_{z \in J(f)} |\psi(z)| + \sup_{z, z' \in J(f), z \neq z'} \frac{|\psi(z) - \psi(z')|}{\dist(z, z')}$ denotes the Lipschitz norm of $\psi$.
Moreover we will say that {\it the \CLT{} holds for} $\nu$, if for every Lipschitz continuous function $\psi : X \to \R$ that is not a coboundary (i.e. it cannot be written in the form $\varphi \circ f - \varphi$) there is $\sigma > 0$ such that for every $x \in \R$ we have,
$$
\nu
\left\{ \frac{1}{\sqrt{n}} \sum_{j = 0}^n \psi \circ f^j < x \right\}
\to
\frac{1}{\sqrt{2\pi\sigma}} \int_{-\infty}^x \exp \left( -\frac{u^2}{2\sigma^2}  \right) du,
\text{ as $n \to \infty$}.
$$
\begin{theoalph}\label{t:statistical properties}
If $f$ is a rational map satisfying the \TCEC, then the invariant measure given by Theorem~\ref{t:invariant measures} is exponentially mixing and the \CLT{} holds for this measure.
\end{theoalph}
To the best of our knowledge this is the first result of this type in the holomorphic setting.
In the case of unimodal maps, a similar result was proved in~\cite{You,KN}.
In~\cite{BLS} this result was shown for Collet-Eckmann multimodal maps having all its critical points of the same critical order.
The proof that we give here for rational maps extends without change for multimodal maps satisfying the \TCEC{} (and hence for those satisfying the \CEC), with no restriction on the critical orders of critical points. 
\subsection{Further results.}
We show that a rational map having an exponentially mixing invariant measure as described in theorems~\ref{t:invariant measures} and~\ref{t:statistical properties}, must satisfy the \TCEC.
This adds yet another equivalent formulation of the \TCEC. 
\begin{theoalph}\label{t:reverse implication}
  Let $f$ be a rational map having an exponentially mixing invariant measure that is absolutely continuous with respect to some conformal measure of~$f$, in such a way that the density is almost everywhere bounded from below by a positive constant.
Then $f$ satisfies the \TCEC.
\end{theoalph}
An analogous result for unimodal maps was shown by T.~Nowicki and D.~Sands in~\cite{NS}.
We follow the same idea of proof, which is to show (by an argument attributed to G.~Keller) that the existence of such a measure implies that the map is uniformly hyperbolic on periodic orbits: there is $\lambda > 1$ such that for every positive integer~$n$ and every repelling periodic point $p$ of $f$ of period~$n$, we have $|(f^n)'(p)| \ge \lambda^n$.
Then we use the result of~\cite{PRS}, that for a complex rational map this last condition is equivalent to the \TCEC.

We also show that the measure given by Theorem~\ref{t:invariant measures} is characterized as the unique invariant measure supported on $J(f)$ whose Hausdorff dimension is equal to $\HD(J(f))$ (Proposition~\ref{p:HD of invariant}) and as the unique equilibrium state with potential $-\HD(J(f)) \ln |f'|$ of~$f$ (Corollary~\ref{c:invariant as equilibrium}).
Similar results for Collet-Eckmann unimodal maps where shown in~\cite{BK}; see also~\cite{PS}.

Finally, we also show that for a rational map satisfying the \TCEC{}, the set of points in the Julia set that are not in the conical Julia set has Hausdorff dimension~$0$.
This is a direct consequence of Theorem~\ref{t:induced map} and Lemma~\ref{l:HD of bad}.
\subsection{Strategy.}
We now explain the strategy of proof of our main results, and simultaneously describe the organization of the paper.

To prove theorems~\ref{t:conformal measures},~\ref{t:invariant measures} and~\ref{t:statistical properties} we use an inducing scheme.
That is, we construct a (Markovian) induced map and then deduce properties of the rational map from properties of the induced map.
To construct a conformal measure supported on the conical Julia set (Theorem~\ref{t:conformal measures}), we basically construct an induced map whose maximal invariant set has the largest possible Hausdorff dimension (equal to the Hausdorff dimension of the Julia set) and which is strongly regular in the sense of~\cite{MUbook}.
Given a rational map having such an induced map, we construct the desired conformal measure in Appendix~\ref{a:conformal via inducing}.
For the existence and statistical properties of the absolutely continuous measure (theorems~\ref{t:invariant measures} and~\ref{t:statistical properties}) we do a ``tail estimate'' and use the results of L.-S.~Young~\cite{You}.
The main difference with previous approaches is that we estimate diameters of pull-backs directly using condition \ExpShrink, and not through derivatives at critical values.

The core of this paper is divided into~$2$ independent parts.
In the first part (\S\S\ref{s:ExpShrink implies DG},\ref{s:nice for TCE}) we construct, for a given rational map satisfying the \TCEC, an induced map which is uniformly hyperbolic in the sense that its derivative is exponentially big with respect to the return time, and that it satisfies some additional properties, see Theorem~\ref{t:induced map} in~\S\ref{ss:constructing nice couples}.
For this, we first show in~\S\ref{s:ExpShrink implies DG} that the \TCEC{} implies a strong form of the ``Backward Contraction'' property of~\cite{Rdecay}.
The construction of the induced map is based on the concept of ``nice couple'' introduced in~\cite{Rdecay}, which is closely related to nice intervals of real one-dimensional dynamics, see~\S\S\ref{ss:nice sets}, \ref{ss:nice couples}.
To each nice couple we associate an induced map of the rational map (\S\ref{ss:canonical induced}).
Then we repeat in~\S\ref{ss:constructing nice couples} the construction of nice couples of~\cite{Rdecay}, using the Backward Contraction property that was shown in~\S\ref{s:ExpShrink implies DG}.
The desired properties of the induced map associated to this nice couple follow easily from the Backward Contraction property.

In the second part (\S\S\ref{s:density}--\ref{s:key lemma}) we give very simple conditions on a nice couple, for a rational map satisfying the \TCEC, so that the associated induced map has the following properties: its maximal invariant set has the largest possible Hausdorff dimension (equal to the Hausdorff dimension of the Julia set) and there is $\alpha \in (0, \HD(J(f)))$ such that,
\begin{equation}\label{e:pressure estimate}
\sum_{W} \diam(W)^\alpha < + \infty,
\end{equation}
where the sum is over all the connected components of the domain of the induced map.
This is stated as the Key Lemma in~\S\ref{s:key lemma}.
Both estimates are very important for the existence of the conformal measure supported on the conical Julia set and for the existence and statistical properties of the absolutely continuous invariant measure.
To prove~\eqref{e:pressure estimate} we introduce in~\S\ref{s:density} a ``discrete density'' that behaves well under univalent and unicritical pull-backs.
Then we prove in~\S\ref{s:nice sets and the density} a general result estimating the discrete density of the domain of the first entry map to a nice set.
The use of the discrete density greatly simplifies the proof of the Key Lemma, as it reduces considerably the combinatorial arguments. 

To prove Theorem~\ref{t:conformal measures} we show a general result in Appendix~\ref{a:conformal via inducing}, that roughly states that if a rational map admits a nice couple satisfying the conclusions of the Key Lemma, then the rational map has a conformal measure supported on the conical Julia set (Theorem~\ref{t:conformal via inducing}).
The proof is very simple: the hypothesis imply that the induced map associated to the nice couple is (strongly) regular in the sense of~\cite{MUbook} and by the results of~\cite{MUbook} we deduce that this induced map has a conformal measure whose dimension is equal to~$\HD(J(f))$.
Then we spread this measure using the rational map, to obtain a conformal measure supported on the conical Julia set.
The proof of Theorem~\ref{t:conformal measures} is in~\S\ref{ss:proof of A}.

The proof of theorems~\ref{t:invariant measures} and~\ref{t:statistical properties} is in~\S\ref{ss:proof of B and C}.
We deduce these results from some results of Young in~\cite{You}, that we briefly recall.
As usual the most difficult part is the ``tail estimate''.
In our case this follows easily from Theorem~\ref{t:induced map} and from the Key Lemma.

In Appendix~\ref{a:induced maps} we gather some general properties of the induced maps considered here, the most important being translations of results of~\cite{MUbook} to our particular setting.

The proof of Theorem~\ref{t:reverse implication} is in~\S\ref{ss:proof of D} and the results about Hausdorff dimension of invariant measures and equilibrium states are in~\S\ref{ss:characterizations}.

\subsection{Notes and references.}\label{si:notes and references}
There are quadratic polynomials having no conformal measures supported on the conical Julia set~\cite{Usurvey}.
Even when there exists such a measure, little is known in general about its properties.
For example, to the best of our knowledge it is not known if the Hausdorff dimension of such a measure is positive.

In~\cite{H} it is shown that for every rational map of degree at least~$2$, the measure of maximal entropy is exponentially mixing and that the \CLT{} is satisfied for this measure.
The last result was first shown in~\cite{DPU}.
However the measure of maximal entropy does not describe well, in general, the geometry of the Julia set: The Hausdorff dimension of the measure of maximal entropy is usually strictly smaller than the Hausdorff dimension of the Julia set~\cite{Zd}.
\subsection{Acknowledgments.}
We are grateful to M.~Urba\'nski for useful conversations and references and to S.~Gou\"ezel for several precisions and references.
The second named author is grateful to the Mathematical Institute of the Polish Academy of Sciences (IMPAN) for hospitality during the preparation of this article.

\section{Preliminaries.}\label{s:preliminaries}
For basic references on the iteration of rational maps, see~\cite{CG, Mi}.

An open subset of the Riemann sphere $\CC$ will be called {\it simply-connected} if it is connected and if its boundary is connected.
\subsection{Spherical metric.}
We will identify the Riemann sphere $\CC$ with $\C \cup \{ \infty \}$ and endow $\CC$ with the spherical metric.
The spherical metric will be normalized in such a way that its density with respect to the Euclidean metric on $\C$ is equal to $z \mapsto (1 + |z|^2)^{-1}$.
With this normalization the diameter of $\CC$ is equal to $\pi /2$.

Distances, balls, diameters and derivatives are all taken with respect to the spherical metric.
For $z \in \CC$ and $r > 0$, we denote by $B(z, r) \subset \CC$ the ball centered at $z$ and with radius $r$.
Note that an open ball of radius $r > \tfrac{1}{2} \diam(\CC)$ is equal to $\CC$.
\subsection{Critical points.}
Fix a complex rational map $f$.
We denote by $\Crit(f)$ the set of critical points of $f$ and by $J(f)$ the Julia set of $f$.
Moreover we put $\sCJ(f) \= \CJ$.
When there is no danger of confusion we denote $\Crit(f)$ and $\sCJ(f)$ just by $\Crit$ and $\sCJ$.

For simplicity we will assume that no critical point in $\sCJ$ is mapped to another critical point under forward iteration.
The general case can be handled by treating whole blocks of critical points as a single critical point; that is, if the critical points $c_0, \ldots, c_k \in J(f)$ are such that $c_i$ is mapped to $c_{i + 1}$ by forward iteration, and maximal with this property, then we treat this block of critical points as a single critical point. 

\subsection{Pull-backs.}
Given a subset $V$ of $\CC$ and an integer $n \ge 0$, the connected components of $f^{-n}(V)$ will be called {\it pull-backs of $V$ by $f^n$}.
Note that the set $V$ is not assumed to be connected.

When considering a pull-back of a ball, we will implicitly assume that this ball is disjoint from the forward orbits of critical points not in $J(f)$.
So such a pull-back can only contain critical points in $J(f)$.
If $f$ does not have indifferent cycles (e.g. if $f$ satisfies the \TCEC), it follows by the Fatou-Sullivan classification of the connected components of the Fatou set~\cite{Bea,CG,Mi}, that there is a neighborhood of $J(f)$ disjoint from the forward orbit of critical points not in $J(f)$.
So in this case every ball centered at a point of $J(f)$ and of sufficiently small radius, will meet our requirement.
\subsection{Distortion of univalent pull-backs.}
We will now state a version of Koebe Distortion Theorem, taking into account that derivatives are taken with respect to the spherical metric.
Given $x \in J(f)$, $r > 0$ and an integer $n \ge 1$, let $W$ be a connected component of $f^{-n}(B(x, r))$ on which $f^n$ is univalent.
Then for $\e \in (0, 1)$ let $W(\e) \subset W$ be the preimage of $B(x, \e r)$ by $f^n$ in~$W$.

Fix two periodic orbits $O_1$ and $O_2$ of $f$ of period at least~$2$ and 
let $r_K > 0$ be sufficiently small such that for every $x \in \CC$ the ball $B(x,r_K)$ is 
disjoint from either $O_1$ or $O_2$.  
So, for every positive integer $n$ and every component $W$ of 
$f^{-n}(B(x,r_K))$, we have,
$$
\diam(\CC \setminus W) \ge 
\min\{\diam (O_1), \diam (O_2)\} > 0.
$$ 
Hence the following version of Koebe Distortion Theorem holds (see also Lemma~$1.2$ of~\cite{P1}).
\begin{generic}[Koebe]
For every $\varepsilon \in (0, 1)$ there exists a 
constant $K(\varepsilon) > 1$ such that, if $x \in J(f)$, $r > 0$, $n \ge 1$, $W$ and $W(\e)$ are as above with $r \in (0, r_K)$, then the distortion of $f^n$ on $W(\e)$ is bounded by $K(\e)$. 
That is, for every $z_1$, $z_2 \in W(\e)$ we have,
$$
|(f^n)'(z_1)|/|(f^n)'(z_2)|\le K(\varepsilon).
$$
Moreover $K(\varepsilon) \rightarrow 1$ as $\varepsilon \rightarrow 0$.
\end{generic}

\section{TCE implies Backward Contraction.}\label{s:ExpShrink implies DG}
The purpose of this section is to prove the following property of rational maps satisfying \ExpShrink{}, which is closely related to the Backward Contraction property of~\cite{Rdecay}.
\begin{prop}\label{p:decay of geometry}
Let $f$ be a rational map satisfying \ExpShrink{} with constant $\la_{\Exp} > 1$.
Then for every $\la \in (1, \la_{\Exp})$ there are constants $\theta \in (0, 1)$ and $\alpha > 0$, such that for every $\delta > 0$ small and every $c \in \sCJ$ there is a constant $\delta(c) \in [\delta, \delta^\theta]$ satisfying the following property.

For every $c$, $c' \in \sCJ$, every integer $n \ge 1$ and every pull-back $W$ of $B(c, \delta^{-\alpha} \delta(c))$ by $f^n$,
$$
\dist(W, c') \le \delta(c') \, \text{ implies } \,
\diam(W) \le \la^{-n} \delta(c').
$$
\end{prop}
After some preliminary lemmas in~\S\S\ref{ss:distortionBD}, \ref{ss:diameter of pull-backs} and~\ref{ss:pull-backs and critical points}, the proof of this proposition is given in~\S\ref{ss:proof of decay of geometry}.
\subsection{Distortion lemma for bounded degree maps.}\label{ss:distortionBD}
The following is a well-known gene\-ral lemma, that is needed in the proof of Proposition~\ref{p:decay of geometry}.
\begin{lemm}\label{l:distortionBD}
Fix $\varepsilon > 0$ and $\hat{r} > 0$ small.
Let $\hW$ be a simply-connected subset of $\CC$ such that $\diam(\CC \setminus \hW) > \varepsilon$.
Let $z \in \CC$, $r \in (0, \hat{r})$, $\varphi : \hW \to B(z, r)$ a holomorphic and proper map and let $W$ be a connected component of $\varphi^{-1}(B(z, r/2))$. 
Then there is a constant $K > 1$, depending only on $\varepsilon$ and on the degree of $\varphi$, such that for every connected set $U \subset B(z, r/2)$ and every connected component $V \subset W$ of $\varphi^{-1}(U)$, we have
$$
\frac{\diam(V)}{\diam(W)} \ge K \frac{\diam(U)}{r}.
$$
\end{lemm}
\begin{proof}
For a set $D \subset \CC$  conformally equivalent to a disc and for $\eta \subset D$ we denote by $\diam_D(\eta)$ the diameter of $\eta$ with respect to the hyperbolic metric of $D$.

Since $r \in (0, \hat{r})$ there is $K_0 > 0$ such that
\begin{equation}\label{e:distortionBD}
\diam(U)/ r \le K_0 \diam_{B(z, r)}(U).
\end{equation}
Moreover the modulus of the annulus $B(z, r) \setminus \ov{B(z, r/2)}$ is bounded from below in terms of $\hat{r}$ only.
Thus the modulus of the annulus $\hW \setminus \ov{W}$ is bounded in terms of $\hat{r}$ and of the degree of $\varphi$ only.
Let $w \in W$ be such that $\varphi(w) = z$ and consider a conformal representation $\psi : \D \to \hW$ such that $\psi(0) = w$.
Let $W' \= \psi^{-1}(W)$ and $V' \= \psi^{-1}(V)$.

Since $\diam(\CC \setminus \hW) > \varepsilon$, by Koebe Distortion Theorem the distortion of $\psi$ on $W'$ (where $\D \subset \C$ is endowed with the Euclidean metric) is bounded by some constant $\kappa_0 > 1$.
Therefore we have,
$$
\diam(V)/\diam(W) \ge \kappa_0 \diam(V')/\diam(W')
\ge \kappa_0 \diam(V').
$$
Since the modulus of $\D \setminus \ov{W'}$ is equal to that of $\hW \setminus \ov{W}$, which is bounded in terms of the degree of $\varphi$ only, there is a constant $\kappa_1 > 0$ such that
$$
\diam(V') \ge \kappa_1 \diam_\D (V') \ge \kappa_1 \diam_{B(z, r)}(U),
$$
where the last inequality follows from Schwarz' Lemma.
So by~\eqref{e:distortionBD} we have,
$$
\frac{\diam(V)}{\diam(W)}
\ge
\kappa_0 \kappa_1 \diam_{B(z, r)}(U)
\ge
K_0^{-1} \kappa_0 \kappa_1 \frac{\diam(U)}{r}.
$$
\end{proof}
\subsection{Diameter of pull-backs.}\label{ss:diameter of pull-backs}
The following is analogous to Lemma~$1.9$ of~\cite{PrHolder}.
\begin{lemm}\label{l:expansion}
Let $f$ be a rational map satisfying \ExpShrink{} with constants $\lambda_{\Exp} > 1$ and $r_0 > 0$.
Then the following assertions hold.
\begin{enumerate}
\item[1.]
There are constants $C_0 > 0$ and $\theta_0 \in (0, 1)$ such that for every $r \in (0, r_0)$, every integer $n \ge 1$, every $x \in J(f)$ and every connected component $W$ of $f^{-n}(B(x, r))$, we have
$$
\diam(W) \le C_0 \lambda_{\Exp}^{-n} r^{\theta_0}.
$$
\item[2.]
For each $\lambda \in (1, \lambda_{\Exp})$ and $\beta \ge 0$ there is a constant $A_0 \= A_0(\la, \beta) > 0$, such that if in addition $n \ge A_0 \ln r^{-1}$, then
$$
\diam(W) \le \lambda^{-n} r^{1 + \beta}.
$$
\end{enumerate}
\end{lemm}
\begin{proof}

\partn{1}
Put $M = \sup |f'|$.
Given $r \in (0, r_0)$ and $x \in J(f)$ let $m \ge 0$ be the integer such that $M^{-(m - 1)} \le r < M^{-m}$ and let $W_0$ be the connected component of $f^{-m}(B(f^m(x), r_0))$ that contains $x$.
Thus,
$$
B(x, r) \subset B(x, M^{-m}) \subset W_0.
$$

Let $n \ge 1$ be an integer, let $W$ be a pull-back of $B(x, r)$ by $f^n$ and let $W_n$ be the corresponding pull-back of $W_0$.
Then, letting $\theta_0 \= \frac{\ln M}{\ln \la_{\Exp}}$ we have
$$
\diam(W) \le \diam(W_n) \le \la_{\Exp}^{-(n + m)}
=
\la_{\Exp}^{-n} M^{- \theta_0 m} 
\le
M^{\theta_0} \la_{\Exp}^{-n} r^{\theta_0}.
$$

\partn{2}
Let $A_0$ be a constant satisfying $A_0 > (1 + \beta - \theta_0) \left( \ln \frac{\la_{\Exp}}{\la} \right)^{-1}$.
So if $n \ge A_0 \ln \frac{1}{r}$ and $r > 0$ is small enough, we have
$$
M^{\theta_0} \lambda_{\Exp}^{-n} r^{\theta_0} 
\le \lambda^{-n}r^{1 + \beta}.
$$
\end{proof}
\subsection{Pull-backs and critical points.}\label{ss:pull-backs and critical points}
Fix a complex rational map $f$ of degree at least~$2$.
Recall that we assume that no critical point in $\sCJ \= \CJ$ is mapped to a critical point under forward iteration (cf. Preliminaries).

The following is Lemma~$1$ of~\cite{PrLyapunov}.
\begin{lemm}\label{l:recurrence}
Put $M = \sup_{\CC} |f'|$.
Then there is a constant $\kappa > 0$ such that for every $c \in \sCJ$ and every $n \ge 1$, we have
$$
\dist(c, f^{-n}(c)) \ge \kappa M^{-n}.
$$
\end{lemm}
\begin{proof}
For a critical point $c$ of $f$ denote by $d(c)$ the local degree of $f$ at $c$.
Let $C \ge 0$ be such that for $\delta > 0$ small and $c \in \sCJ$, we have $\diam(f(B(c, \delta))) \le C \delta^{d(c)}$.
Given $n \ge 1$ let $w \in f^{-n}(c)$ be a point closest to $c$ and put $r \= \dist(c, w)$.
Letting $C_0(c) \= M^{-1}C3^{d(c)}$ and considering that $B(w, 2r) \subset B(c, 3r)$, we have
$$
\diam(f^n(B(w, 2r))) \le M^{n - 1} \diam(f(B(w,2r)))
\le C_0(c) M^n r^{d(c)}.
$$
Thus $f^n(B(c, r)) \subset f^n(B(w, 2r)) \subset B(c, C_0(c) M^n r^{d(c)})$.
Since $c \in J(f)$ we must have $C_0(c) M^n r^{d(c)} > r$, so the lemma follows with constant $\kappa \= \max_{c \in \sCJ} C_0(c)^{-1}$.
\end{proof}

\

\begin{lemm}\label{l:recurrenceUHP}
Let $f$ be a rational map satisfying \ExpShrink.
Then there are constants $r_1 > 0$ and $A_1 > 0$ such that for every $c \in \sCJ$, every $r \in (0, r_1)$, every $n \ge 1$ and every connected component $W$ of $f^{-n}(B(c, r))$,
$$
\dist(c, W) \le r \, \text{ implies } \, n \ge A_0 \ln 1/r.
$$
\end{lemm}
\begin{proof}
Let $C_0$, $r_0 > 0$ and $\theta_0 \in (0, 1)$ be given by Lemma~\ref{l:expansion}.
So, assuming $r \in (0, r_0)$ we have $\diam(W) \le C_0 r^{\theta_0}$.
Since $W$ contains a $n$-th preimage of $c$ we have by Lemma~\ref{l:recurrence},
$$
\kappa M^{-n} \le r + \diam(W) \le 2C_0 r^{\theta_0}.
$$
Hence the Lemma follows with any constant $A_1 \in (0, \theta_0/\ln M)$, by taking $r_1$ sufficiently small.
\end{proof}
\begin{lemm}\label{l:boundeddegree}
Let $f$ be a rational map satisfying \ExpShrink{} and let $N \ge 1$ be the number of critical points of $f$ in $J(f)$.
Then there are constants $r_1 > 0$, $\xi_1 > 1$ and $\nu_1 > 0$ such that for every $x \in J(f)$, every $r \in (0, r_1)$ and every pull-back $W_0 \= B(x, r)$, $W_1, \ldots$ the following assertions hold.
\begin{enumerate}
\item[1.]
The sequence $0 \le m_0 < m_1 < \ldots$ of all integers such that $W_{m_i}$ contains a critical point satisfies,
$$
m_i \ge \nu_1 \xi_1^i \ln 1/r, \, \text{ for } i > N.
$$
\item[2.]
For every $A > 0$ there is a constant $D_1 \= D_1(A) \ge 1$ such that if $n \le A \ln 1/r$ then the degree of 
$$
f^n : W \to B(x, r)
$$
is at most $D_1$.
\end{enumerate}
\end{lemm}
\begin{proof}
Note that part~$2$ is a direct consequence of part~$1$.

Fix $\lambda \in (1, \la_{\Exp})$ and let $r_0 > 0$, $C_0 > 0$ and $\theta_0 \in (0, 1)$ be as in Lemma~\ref{l:expansion}, so $\diam(W_i) \le C_0 \la^{-i} r^{\theta_0}$ for $i \ge 0$.
In particular we may suppose that $W_{m_i}$ contains a unique critical point $c_i \in \sCJ$.

Fix $i > N$.
Then there is a critical point $c$ such that
$$
k = \# \{ j \in \{ 0, \ldots, i \} \mid c_j = c \} \ge i/N > 1;
$$
let $m_0 \le n_0  < n_1 < \ldots < n_{k - 1} \le m_i$ be all the integers such that $c \in W_{n_j}$.
Thus for every $j = 1, \ldots, k - 1$ the set $W_{n_j}$ contains a $(n_j - n_{j-1})$-th preimage of $c$.
Then Lemma~\ref{l:recurrence} implies that,
$$
\kappa M^{-(n_j - n_{j - 1})} \le \diam(W_{n_j}) 
\le C_0\la^{-n_j}r^{\theta_0}.
$$
So letting $\xi \= \frac{\ln M}{\ln(M/\la)} > 1$ and $\nu \= \tfrac{1}{2} \frac{\theta_0}{\ln(M/\la)}$ and assuming $r_1 > 0$ small enough, we have $n_j \ge \xi n_{j - 1} + \nu \ln 1/r$.
Therefore,
$$
m_i \ge n_k \ge \nu \xi^{k - 2} \ln 1/r
\ge (\nu \xi^{-2}) (\xi^{1/N})^i \ln 1/r,
$$
and the lemma follows with $\xi_1 \= \xi^{1/N}$ and $\nu_1 \= \nu \xi^{-2}$.
\end{proof}
\subsection{Proof of Proposition~\ref{p:decay of geometry}.}\label{ss:proof of decay of geometry}
Let $f$ be a rational map satisfying condition \ExpShrink{} with constants $\lambda_{\Exp} > 1$ and $r_0 > 0$ and let $C_0 > 0$ and $\theta_0 \in (0, 1)$ be the constants given by part~$1$ of Lemma~\ref{l:expansion}.
Moreover let $r_1 > 0$ and $A_1 > 0$ be the constants given by Lemma~\ref{l:recurrenceUHP}.
Denote by $N$ the cardinality of $\sCJ$ and choose $\la_0 \in (1, \la_{\Exp})$,
$$
\alpha_2 \in \left( 0, \min \left\{ (A_1 \theta_2 / 4N) \ln (\la_0/\la), \theta_0/2N \right\} \right),
$$
$$
\theta_2  \in (0, \theta_0 - 2N\alpha_2).
$$

Choose $\la \in (1, \la_0)$ and fix $\delta > 0$ small enough so that $\delta^{\theta_2} < r_1$.
Let $c_0 \in \sCJ$ and consider successive pull-backs $W_0 \= B(c_0, \delta), W_1, \ldots$ \,.
For $j \ge 0$ define $c_j \in \sCJ$, $\delta_j > 0$ and $m_j \ge 1$ inductively as follows.
Put $\delta_0 \= \delta$ and $m_0 \= 0$.
Once $c_{j - 1}$, $\delta_{j - 1}$ and $m_{j - 1}$ are defined, let $m_j > m_{j - 1}$ be the least integer such that for some $c_j \in \sCJ$ we have
$$
\delta_j \= \la^{m_j}\delta^{- j\alpha_2} \diam(W_{m_j})
> \dist(c_j, W_{m_j}).
$$
\begin{lemm}\label{l:uniform}
There is a constant $A_2 > 0$ independent of $\delta$, such that the following assertions hold.
\begin{enumerate}
\item[1.]
For every $j \ge 0$ we have $\delta_j \le \delta^{\theta_2}$.
\item[2.]
If $j \ge 1$ is such that $m_j \ge A_2 \ln 1/\delta$, then $\delta_j \le \delta$.

\end{enumerate}
\end{lemm}
\begin{proof}
Part~$1$ and~$2$ of the lemma will be shown in parts~3 and~4 below.

\partn{1}
If $\delta > 0$ is sufficiently small, then for every $j = 0, \ldots, 2N$ we have
\begin{multline*}
\delta_j
=
\la^{m_j} \delta^{-j\alpha_2} \diam(W_{m_j})
\le
\la^{m_j} \delta^{-j\alpha_2} C_0 \la_0^{-m_j}\delta^{\theta_0}
\\ \le
C_0 \delta^{\theta_0 - 2N \alpha_2}
<
\delta^{\theta_2}/2.
\end{multline*}

\partn{2}
Consider $j' > j \ge 0$ such that $c_{j'} = c_{j}$ and such that
$$
\delta' \= \max \{ \delta_{j}, \delta_{j'} \} 
\le \delta^{\theta_2}/2 < r_1/2.
$$
Then $W_{m_j}, W_{m_{j'}} \subset B(c_j, 2\delta')$ and the pull-back $W$ of $B(c_j, 2\delta')$ by $f^{m_{j'} - m_j}$ con\-taining $W_{m_{j'}}$ intersects $B(c_j, 2\delta')$.
So by Lemma~\ref{l:recurrenceUHP}
\begin{equation*}
m_{j'} - m_{j} \ge A_1 \ln 1/(2\delta') \ge A_1 \theta_2 \ln 1/\delta.
\end{equation*}

\partn{3}
Given $k \ge 2N$ suppose by induction that $\delta_j \le \delta^{\theta_2}$ for $j = 0, \ldots, k$. 
By part~$1$ this holds for $k = 2N$.
Notice that there must be a critical point $c \in \sCJ$ such that,
$$
\# \{ j \in \{1, \ldots, k \} \mid c_j = c \} \ge k/N.
$$
Since $k \ge 2N$ it follows by part~$2$ that 
\begin{equation}\label{e:time versus hits}
m_k
\ge
(k/N - 1) A_1 \theta_2 \ln 1/\delta
\ge
(A_1 \theta_2 / 2N) k \ln 1/\delta.
\end{equation}
By definition $\alpha_2 < (A_1 \theta_2/4N) \ln (\la_0/\la)$, so
$
(\la_0/\la)^{m_k}
>
\delta^{- 2k \alpha_2}.
$
Considering that $m_{k + 1} > m_k$ we have
\begin{multline}
\label{e:uniform1}
\diam(W_{m_{k + 1}})
\le
C_0 \la_0^{-m_{k + 1}} \delta^{\theta_0}
\le \\ \le 
C_0 \la^{-m_{k + 1}} \delta^{2k \alpha_2} \delta^{\theta_0}
<
\la^{-m_{k+1}} \delta^{(k + 1)\alpha_2} \delta^{\theta_0},
\end{multline}
therefore
$$
\delta_{k + 1} = \la^{m_{k + 1}} \delta^{-(k + 1) \alpha_2} \diam(W_{m_{k + 1}})
\le \delta^{\theta_0} < \delta^{\theta_2}.
$$

\partn{4}
To prove part~$2$ of the lemma, let $A_0(\la_0, 0)$ be the constant given by part~$2$ of Lemma~\ref{l:expansion} with $\la \= \la_0$ and $\beta \= 0$ and put
$$
A_2 \= \max \{ A_0(\la_0, 0), A_1 \theta_2 \}.
$$
Let $j \ge 1$ be such that $m_j \ge A_2 \ln 1/\delta$.
If $k \= j - 1 \ge 2N$, then by Lemma~\ref{l:expansion} and by the same reasoning as in part~3, we conclude that inequality \eqref{e:uniform1} holds with $\delta^{\theta_0}$ replaced by $\delta$.
So in this case we have $\delta_j = \delta_{k + 1} \le \delta$.
If $k \= j - 1 < 2N$, then by definition of $A_2$ we have $m_k \ge A_1 \theta_2 \ln 1/\delta \ge (A_1 \theta_2 / 2N) k \ln 1/\delta$.
So again inequality~\eqref{e:uniform1} holds with $\delta^{\theta_0}$ replaced by $\delta$ and therefore $\delta_j = \delta_{k + 1} \le \delta$.
\end{proof}

\begin{proof}[Proof of Proposition~\ref{p:decay of geometry}]
For $c \in \sCJ$ put
$$
\delta(c) \= \sup_{ W_0, W_1, \ldots, \,\, c_j = c}
\la^{m_j}\delta^{-j\alpha_2} \diam(W_j),
$$
where the supremum is taken over all $c_0 \in \sCJ$ and all successive pull-backs $W_0 \= B(c_0, \delta), W_1, \ldots$, where $c_j$, $\delta_j$ and $m_j$ are defined as above.

Put  $\theta \= \theta_2$ and fix $\alpha \in (0, \alpha_2)$.
In what follows we will prove that, if $\delta > 0$ is sufficiently small, then the assertion of Proposition~\ref{p:decay of geometry} holds for this choice $\theta$, $\alpha$ and $\delta(c)$.

\partn{1}
Taking $c_0 = c$ and $j = 0$ in the definition of $\delta(c)$, we have $\delta(c) \ge \delta$.
Moreover by part~$1$ of Lemma~\ref{l:uniform} we have $\delta(c) \le \delta^{\theta_2}$.
On the other hand by part~$2$ of the same lemma, $\delta_j > \delta$ implies that $m_j \le A_2 \ln 1/\delta$.
So the supremum defining $\delta(c)$ is realized.

\partn{2}
Let $c$, $c' \in \sCJ$, $n \ge 1$ and let $W$ be a pull-back of $B(c, \delta^{-\alpha} \delta(c))$ by $f^n$ such that $\dist(c', W) \le \delta(c')$.
By part~$2$ of Lemma~\ref{l:expansion} it follows that there is $A_0 > 0$ such that, if $n \ge A_0 \ln 1/\delta$, then $\diam(W) \le \la^{-n} \delta \le \la^{-n} \delta(c')$.
So we assume that $n \le A_0 \ln 1/\delta$.

By part~$1$ there is $c_0 \in \sCJ$ and a pull-back $W_0 \= B(c_0, \delta), W_1, \ldots, W_{m_k}$, such that for some $k \ge 0$ we have $c_k = c$ and
$$
\delta(c)
=
\delta_k 
= \la^{m_k} \delta^{- k\alpha_2} \diam(W_{m_k})
> \dist(c, W_{m_k}).
$$
Note that $W_{m_k} \subset B(c, 2\delta(c)) \subset B(c, \delta^{-\alpha} \delta(c))$.
Consider the successive pull-backs $W_{m_k}, W_{m_k + 1}, \ldots, W_{m_k + n}$, such that $W_{m_k + n} \subset W$.
Let $c_j$, $\delta_j$ and $m_j$ be as in Lemma~\ref{l:uniform} for the pull-back $W_0, \ldots, W_{m_k + n}$.
Then for some $\ell > k$ we have $c_\ell = c'$ and $m_\ell = m_k + n$.
Hence by definition of $\delta(c')$,
\begin{equation}\label{e:uniform2}
\diam(W_{m_k + n}) \le \la^{-(m_k + n)} \delta^{\ell \alpha_2} \delta(c').
\end{equation}

\partn{3} 
Let $\hW$ be the pull-back of $B(c, 2 \delta^{-\alpha} \delta(c))$ by $f^n$ that contains $W$.
Since $n \le A_0 \ln 1/\delta$ it follows by Lemma~\ref{l:boundeddegree} that the degree of $f^n : \hW \to B(c, 2 \delta^{-\alpha} \delta(c))$ is bounded by a constant depending on $A_0$ only.
Hence by Lemma~\ref{l:distortionBD} there is a constant $K > 1$ such that,
$$
\diam(W)/\diam(W_{m_k + n}) \le K \delta^{-\alpha}\delta(c)/\diam(W_{m_k})
= K \delta^{-\alpha} \la^{m_k} \delta^{- k \alpha_2}.
$$
So by ($\ref{e:uniform2}$) and assuming $\delta > 0$ small enough we have,
$$
\diam(W) \le K \delta^{-\alpha} \la^{-n}\delta^{(\ell - k)\alpha_2} \delta(c')
< \la^{-n} \delta(c').
$$

This ends the proof of Proposition~\ref{p:decay of geometry}.
\end{proof}

\section{Induced maps.}\label{s:nice for TCE}
Given a rational map satisfying the \TCEC, the purpose of this section is to construct an induced map that is hyperbolic in the sense that its derivative is exponentially big with respect to the return time, and that it satisfies some additional properties (Theorem~\ref{t:induced map}).
The construction of this induced map is based on the construction of ``nice couples'' in~\cite{Rdecay}.
We first recall the definition of nice sets (\S\ref{ss:nice sets}) and of nice couples (\S\ref{ss:nice couples}), and then we explain how to associate to each nice couple an induced map (\S\ref{ss:canonical induced}).
The statement and proof of Theorem~\ref{t:induced map} is in~\S\ref{ss:constructing nice couples}.
\subsection{Nice sets.}\label{ss:nice sets}
Let $f$ be a complex rational map.
We will say that a neighborhood $V$ of $\sCJ$ that is disjoint from the forward orbits of critical points not in $\sCJ$ is a {\it nice set for} $f$, if it satisfies the following properties.
The set $V$ is the union of sets $V^c$, for $c \in \sCJ$, such that $V^c$ is a simply-connected neighborhood of $c$, such that the closures of the sets $V^c$ are pairwise disjoint and such that for every pull-back $W$ of $V$ we have either
$$
\ov{W} \cap \ov{V} = \emptyset \, \text{ or } \, \ov{W} \subset V.
$$
Let $V = \cup_{c \in \sCJ} V^c$ be a nice set for $f$.
Note that if $W$ and $W'$ are distinct pull-backs of $V$, then we have either,
$$
\ov{W} \cap \ov{W'} = \emptyset,
\, \ov{W} \subset W'
\, \text{ or } \,
\ov{W'} \subset W.
$$
For a pull-back $W$ of $V$ we denote by $c(W)$ the critical point in $\sCJ$ and by $m_W \ge 0$ the integer such that $f^{m_W}(W) = V^{c(W)}$.
Moreover we put,
$$
K(V)
=
\{ z \in \CC \mid \text{ for every $n \ge 0$ we have $f^n(z) \not \in V$} \}.
$$
Note that $K(V)$ is a compact and forward invariant set and for each $c \in \sCJ$ the set $V^c$ is a connected component of $\CC \setminus K(V)$.
Moreover, if $W$ is a connected component of $\CC \setminus K(V)$ different from the $V^c$, then $f(W)$ is again a connected component of $\CC \setminus K(V)$.
It follows that $W$ is a pull-back of $V$ and that $f^{m_W}$ is univalent on $W$.
\subsection{Nice couples.}\label{ss:nice couples}
A {\it nice couple for} $f$ is a couple $(\hV, V)$ of nice sets for $f$ such that for every pull-back $\hW$ of $\hV$ we have either
$$
\ov{\hW} \cap \ov{V} = \emptyset
\, \text{ or } \,
\ov{\hW} \subset V.
$$
Let $(\hV, V)$ be a nice couple for $f$.
Then for each pull-back $W$ of $V$ we denote by $\hW$ the corresponding pull-back of $\hV$, in such a way that $W \subset \hW$, $m_{\hW} = m_W$ and $c(\hW) = c(W)$.
If $W$ and $W'$ are disjoint pull-backs of $V$ such that the sets $\hW$ and $\widehat{W'}$ intersect, then we have either
\begin{equation}\label{e:Markov property for pull-backs}
\hW \subset \widehat{W'} \setminus W'
\, \text{ or } \,
\widehat{W'} \subset \hW \setminus W.
\end{equation}
If $W$ is a connected component of $\CC \setminus K(V)$, then for every $j = 0, \ldots, m_W - 1$, the set $f^j(W)$ is a connected component of $\CC \setminus K(V)$ different from the $V^c$, and $\widehat{f^j(W)}$ is disjoint from $V$.
It follows that $\widehat{f^j(W)}$ does not contain critical points of $f$ and that $f^{m_W}$ is univalent on $\hW$.

\subsection{The canonical induced map associated to a nice couple.}\label{ss:canonical induced}
Let $f$ be a complex rational map and let $(\hV, V)$ be a nice couple for~$f$.
We will say that an integer $m \ge 1$ is a {\it good time for a point $z$} in $V$, if $f^m(z) \in V$ and if the pull-back of $\hV$ by $f^m$ to $z$ is univalent.
Let $D$ be the set of all those points in $V$ having a good time and for $z \in D$ denote by $m(z) \ge 1$ the least good time of $z$.
Then the map $F : D \to V$ defined by $F(z) \= f^{m(z)}(z)$ is called {\it the canonical induced map associated to $(\hV, V)$}.
We denote by $J(F)$ the maximal invariant set of~$F$.

As $V$ is a nice set, it follows that each connected component $W$ of~$D$ is a pull-back of $V$.
Moreover, $f^{m_W}$ is univalent on $\hW$ and for each $z \in W$ we have $m(z) = m_W$.
Similarly, for each positive integer~$n$, each connected component $W$ of the domain of definition of $F^n$ is a pull-back of $V$ and $f^{m_W}$ is univalent on $\hW$.
Conversely, if~$W$ is a pull-back of $V$ contained in $V$ and such that $f^{m_W}$ is univalent on $\hW$, then there is $c \in \sC$ and a positive integer~$n$ such that $F^n$ is defined on~$W$ and $F^n(W) = V^c$.
In fact, in this case $m_W$ is a good time for each element of $W$ and therefore $W \subset D$.
Thus, either we have $F(W) = V^{c(W)}$, and then $W$ is a connected component of $D$, or $F(W)$ is a pull-back of $V$ contained in $V$ such that $f^{m_{F(W)}}$ is univalent on $\widehat{F(W)}$.
Thus, repeating this argument we can show by induction that there is a positive integer~$n$ such that $W$ is defined on $W$ and that $F^n(W) = V^{c(W)}$.
\begin{lemm}\label{l:mixingness}
For every rational map $f$ there is $r > 0$ such that if $(\hV, V)$ is a nice couple satisfying
\begin{equation}\label{e:diam bound}
\max_{c \in \sC} \diam(\hV^c) \le r,
\end{equation}
then the canonical induced map $F : D \to V$ associated to $(\hV, V)$ is topologically mixing on $J(F)$.
Moreover there is $\tc \in \sC$ such that the set
\begin{equation}\label{e:return times}
\{ m_W \mid W \text{ c.c. of $D$ contained in $V^{\tc}$ such that $F(W) = V^{\tc}$} \}
\end{equation}
is non-empty and its greatest common divisor is equal to~$1$.
\end{lemm}
\begin{proof}
Let $p$ be a repelling periodic point of~$f$.
By the locally eventually onto property of Julia sets~\cite{CG, Mi}, for each critical point $c \in \sC$ there is a backward orbit starting at~$c$ and that is asymptotic to the backward periodic orbit of~$f$ starting at~$p$.
As our standing assumption is that no critical point is mapped into another critical point under forward iteration, this backward orbit does not contains critical points.
Let $r > 0$ be sufficiently small, so that $B(\sC, r)$ intersects each of these backward orbits only at its starting point in~$\sC$.

Let $(\hV, V)$ be a nice couple for~$f$ satisfying~\eqref{e:diam bound} and let $F$ be the canonical induced map associated to $(\hV, V)$.
Since the partition of $J(F)$ induced by the connected components of~$D$ is generating, it follows that for every open set~$U$ intersecting $J(F)$ there is $c \in \sC$ and an integer $n$ such that $F^n(U)$ contains $V^c$.
Thus, to prove that $F$ is topologically mixing we just have to show that for every $c, c' \in \sC$ there is $n_0 > 0$ such that for every positive integer $n \ge n_0$ the set $F^n(V^{c'})$ contains~$V^c$.
For this, we will show that for every $c, c' \in \sC$ there is a positive integer~$n$ such that the set $F^n(V^{c'})$ contains $V^c$, and then that there is $c_0 \in \sC$ such that $F(V^{c_0})$ contains $V^{c_0}$.

By~\eqref{e:diam bound} it follows that for each $c \in \sC$, the periodic point~$p$ is accumulated by a sequence $(\hW_n)_{n \ge 1}$ of connected components of~$\CC \setminus K(\hV)$ that are pull-backs of $\hV^c$.
By the locally eventually onto property of Julia sets, it follows that for each $c' \in \sC$ there is an integer $n(c')$ and $q(c') \in V^{c'}$ such that $f^{n(c')}(q(c')) = p$.
As for large~$n$ the set $\hW_n$ is disjoint from $f^{n(c')}(\sC)$, it follows that for large~$n$ the pull-back of $\hW_n$ by $f^{n(c')}$ near $q(c')$ is a univalent pull-back of~$\hV^c$.
This shows that for every $c, c' \in \sC$ there is a positive integer $n$ such that $F^n(V^{c'})$ contains $V^c$.

Let us prove now that there is $c_0 \in \sC$ such that $F(V^{c_0})$ contains $V^{c_0}$.
Choose an arbitrary $c' \in \sC$ and let $q(c')$ as before.
Let $n \ge 0$ be the largest integer such that $q \= f^n(q(c')) \in V$ and let $c_0 \in \sC$ be such that $q \in V^{c_0}$.
Then $q$ is an iterated preimage of~$p$ such that $f(q) \in K(V)$.
As before, for every $c \in \sC$ the point $q$ is accumulated by univalent pull-backs of $\hV^c$.
As $f(q) \in K(V)$, it follows that the corresponding pull-backs of $V^c$ are contained in the domain of~$F$ and thus that $F(V^{c_0})$ contains $V^c$.
In particular $F(V^{c_0})$ contains $V^{c_0}$.

To prove the final statement we will use the fact that every rational map has a repelling periodic point of each sufficiently large period.
This follows from the result of~\cite[Theorem~$6.2.2$, p.102]{Bea}, that every rational map has a periodic point of a given (minimal) period greater or equal than~$4$, and from the fact that a rational map possesses at most finitely many non-repelling periodic points.
For the proof of the final statement, notice that we can take the sequence $(\hW_n)_{n \ge 1}$ above, in such a way that $(m_{\hW_n})_{n \ge 1}$ is an arithmetic progression for which the difference between~$2$ consecutive terms is equal to the minimal period of~$p$.
This shows that the set~\eqref{e:return times} contains a set of the form $\{a + n b \mid n \ge 0 \}$, where $b$ is the minimal period of~$p$.
Repeating the argument with $\#\sC + 1$ repelling periodic points whose periods are pairwise distinct prime numbers, we can find $\tc \in \sC$ for which the set~\eqref{e:return times} contains sets of the form $\{a_0 + n b_0 \mid n \ge 0 \}$ and $\{a_1 + n b_1 \mid n \ge 0\}$, where $b_0$ and $b_1$ are distinct prime numbers.
This implies that the greatest common divisor of the set~\eqref{e:return times} is equal to~$1$.
\end{proof}
\subsection{Constructing nice couples.}\label{ss:constructing nice couples}
The purpose of this section is to prove the following result.
\begin{theoalph}\label{t:induced map}
Let $f$ be a rational map satisfying \ExpShrink{} with constant $\la_{\Exp} > 1$.
Then for every $\la \in (1, \la_{\Exp})$, $\modul > 0$ and $r > 0$ there is a nice couple $(\hV, V)$ such that,
\begin{equation}
\min_{c \in \sCJ} \modulus (\hV^c \setminus \ov{V^c}) \ge \modul,
\max_{c \in \sCJ} \diam(\hV^c) \le r,
\end{equation}
and such that the canonical induced map $F : D \to V$ associated to $(\hV, V)$ satisfies the following property: For every $z \in D$ we have $|F'(z)| \ge \la^{m(z)}$.
\end{theoalph}
In view of Proposition~\ref{p:decay of geometry}, this theorem is a direct consequence of the following proposition.
\begin{prop}\label{p:nice couples}
Let $f$ be a rational map satisfying \ExpShrink{} with constant $\lambda_{\Exp} > 1$ and choose $\lambda \in (1, \lambda_{\Exp})$ and $\tau \in (0, \tfrac{1}{4})$.
For $\delta > 0$ small and $c \in \sCJ$, let $\delta(c) \ge \delta$ be given by Proposition~\ref{p:decay of geometry}.
Then for every $\delta > 0$ sufficiently small there is a nice couple $(\hV, V)$ for $f$ such that for each $c \in \sCJ$ we have,
\begin{equation}\label{e:nice couple}
B(c, \tfrac{1}{2}\delta(c)) \subset \hV^c \subset B(c, \delta(c))
\text{ and }
B(c, \tau \delta(c)) \subset V^c \subset B(c, 2\tau\delta(c)).
\end{equation}
\end{prop}
The proof of this proposition is a repetition of~\cite[Proposition~$6.6$]{Rdecay}.
We include it here for completeness.
It depends on the following lemma.
\begin{lemm}\label{l:prenice}
Given $\delta > 0$ small, put $\tV_0 = \cup_{c \in \sCJ} B(c, \delta(c))$ and
$$
K(\tV_0)
= 
\{ z \in \CC \mid \text{ for every $n \ge 0$ we have $f^n(z) \not \in \tV_0$} \}.
$$
If $\delta > 0$ is sufficiently small, then for each $c \in \sCJ$ the connected component $\tV^c$ of $\CC \setminus K(\tV_0)$ that contains $c$ satisfies,
$$
B(c, \delta(c)) \subset \tV^c \subset B(c, 2\delta(c)).
$$
\end{lemm}
\begin{proof}
Let $\alpha > 0$ be given by Proposition~\ref{p:decay of geometry} and suppose that $\delta > 0$ is sufficiently small so that $\delta^{-\alpha} > 2$.

Given $c \in \sCJ$ and an integer $n \ge 0$, let $\tV_n^c$ be the connected component of $\cup_{j = 0, \ldots, n - 1} f^{-j}(\tV_0)$ that contains $c$.
Note that $\tV^c_0 = B(c, \delta(c))$ and that $\tV_0 = \cup_{c \in \sCJ} \tV_0^c$.
Moreover, note that $\tV^c_n$ is increasing with $n$ and that $\tV^c = \cup_{n \ge 0} \tV_n^c$.
To prove the lemma is enough to show that for every integer $n \ge 0$ we have $\tV_n^{c} \subset \tB(c, 2\delta(c))$.

We will proceed by induction in~$n$.
The case $n = 0$ being trivial, suppose by induction hypothesis that the assertion holds for some $n \ge 0$ and fix $c \in \sCJ$.
We will show that the assertion holds for $n + 1$.
For every point $z \in \tV_{n + 1}^{c}$ there is an integer $m \in \{ 0, \ldots, n + 1 \}$ and $c_0 \in \sCJ$ such that $f^m(z) \in B(c_0, \delta(c_0))$; let $m(z)$ be the least of such integers.
Let $X$ be a connected component of $\tV_{n + 1}^{c} \setminus B(c, \delta(c))$ and let $z \in X$ for which $m(z)$ is minimal among points in $X$.
Let $c_0 \in \sCJ$ be such that $f^{m(z)}(z) \in B(c_0, \delta(c_0))$.
Considering that $m(z) > 0$, we have by induction hypothesis 
$$
f^{m(z)}(X) \subset \tV_n^{c_0} \subset \tB(c_0, 2\delta(c_0)).
$$
Then Proposition~\ref{p:decay of geometry} implies that $\diam(X) < \delta(c)$.
Thus $X \subset B(c, 2\delta(c))$ and $\tV_{n + 1}^{c} \subset B(c, 2 \delta(c))$.
This completes the induction step and the proof of the lemma.
\end{proof}

\begin{proof}[Proof of Proposition~\ref{p:nice couples}.]
Let $\theta \in (0, 1)$ and $\alpha > 0$ be given by Proposition~\ref{p:decay of geometry} and choose $\delta > 0$ sufficiently small so that the conclusions of Proposition~\ref{p:decay of geometry} hold and so that $\delta^{-\alpha} > 2$.
Furthermore, we assume that $\delta > 0$ is sufficiently small so that the least positive integer $L \ge 1$ such that $f^L(\sCJ)$ intersects $B(\sCJ, \delta^\theta)$, satisfies $\lambda^{-L} < \tau$.

For $c \in \sCJ$ let $\tV^c$ be the connected component of $\CC \setminus K(\tV_0)$ that contains~$c$ and put $\tV = \cup_{c \in \sCJ} \tV^c$.
Then Lemma~\ref{l:prenice} implies that for each $c \in \sCJ$ we have $\tV^c \subset B(c, 2\delta(c)) \subset B(c, \delta^{-\alpha}\delta(c))$.

In part~$1$ below, for each $c \in \sCJ$ we construct sets $\hV^c$ and $V^c$ satisfying~\eqref{e:nice couple} and such that in addition $\ov{\hV^c} \subset \tV^c$ and $\partial \hV^c, \partial V^c \subset f^{-1}(K(\tV_0))$.
In part~$2$ we conclude from these properties that the sets $\hV \= \cup_{c \in \sCJ} \hV^c$ and $V \= \cup_{c \in \sCJ} V^c$ are nice sets for~$f$ and that $(\hV, V)$ is a nice couple for~$f$.

\partn{1}
Note that Proposition~\ref{p:decay of geometry} implies that for every $c \in \sCJ$ and every pull-back $W$ of $\tV$ by $f$ that intersects $B(c, \delta(c))$, we have $\diam(W) < \tau \delta(c)$.

We will just construct $\hV^c$, the construction of $V^c$ is analogous.
For $c \in \sCJ$ put,
\begin{multline*}
\ddot{V}^c = B(c, \tfrac{1}{2} \delta(c)) \cup
\left\{ \text{ connected components of $\CC \setminus f^{-1}(K(\tV_0))$} \right.
\\
\left.  \text{intersecting } B(c, \tfrac{1}{2} \delta(c)) \right\}.
\end{multline*}
We have $B(c, \tfrac{1}{2} \delta(c)) \subset \ddot{V}^c \subset B(c, (\tfrac{1}{2} + \tau) \delta(c))$ and $\partial \ddot{V}^c \subset f^{-1}(K(\tV_0))$.
Now put,
\begin{multline*}
\hV^c = \ddot{V}^c \cup
\left\{ \text{ connected components $W$ of $\CC \setminus \ddot{V}^c$} \right.
\\
\left.  \text{such that } \diam(W) < \diam(\CC) \right\}.
\end{multline*}
It follows that $\hV^c$ is simply-connected, that $B(c, \tfrac{1}{2} \delta(c)) \subset \hV^c \subset B(c, (\tfrac{1}{2} + \tau) \delta(c))$ and that $\partial \hV^c \subset \partial \ddot{V}^c \subset f^{-1}(K(\tV_0))$.
Observe finally that $\ov{\hV^c} \subset \tV^c$.

\partn{2}
Given $c \in \sCJ$ let $W_0$ be equal to either $\hV^c$ or $V^c$ and let $W_0, W_1, \ldots$ be successive pull-backs by $f$.
For an integer $n \ge 1$ let $\tW_n$ be the connected component of $\CC \setminus K(\tV_0)$ that contains $W_n$.
Since $\partial \hV \subset f^{-1}(K(\tV_0))$ (resp. $\partial V \subset f^{-1}(K(\tV_0))$), it follows that we have either $\tW_n \subset \hV$ or $\tW_n \cap \hV = \emptyset$ (resp. $\tW_n \subset V$ or $\tW_n \cap V = \emptyset$).
So to prove that the sets $\hV$ and $V$ are nice sets for $f$ and that $(\hV, V)$ is a nice couple for $f$, is enough to prove that for each $n \ge 0$ we have $\ov{W_n} \subset \tW_n$.

We proceed by induction.
For $n = 0$ just note that $\ov{W_0} \subset \tW_0$ we have, because $\ov{V} \subset \ov{\hV} \subset \tV_0$.
Suppose by induction hypothesis, that for some integer $n \ge 1$ we have $\ov{W_n} \subset \tW_n$.
If $\tW_n$ does not intersect $\sCJ$ then the map $f : \tW_n \to \tW_{n - 1}$ is proper, so we have $\ov{W_n} \subset \tW_n$ in this case by the induction hypothesis.
If $\tW_n$ intersects $\sCJ$ then let $\tW_n'$ be the connected component of $f^{-1}(\tW_{n - 1})$ that contains $W_n$, so that $\tW_n' \subset \tW_n$.
By the induction hypothesis, we have $\ov{W_n} \subset \tW_n' \subset \tW_n$.
This completes the induction step and ends the proof of the lemma.
\end{proof}

\section{The density $\Na{\cdot}$.}\label{s:density}
Throughout all this section we fix $\alpha > 0$.
\subsection{Families of subsets of $\CC$ and the density~$\Na{\cdot}$.}
Given a family $\fF$ of subsets of $\CC$ put
$$
\supp(\fF) \= \cup_{W \in \fF} W \ \text{ and}
$$
$$
\Na{\fF}
\=
\sup_{\varphi} \left( \sum_{W \in \fF} \diam(\varphi(W))^\alpha \right),
$$
where the supremum is taken over all M\"obius transformations $\varphi$.
For two such families $\fF$ and $\fF'$ we have,
$$
\Na{\fF \cup \fF'} \le \Na{\fF} + \Na{\fF'}.
$$
In the case $\fF \subset \fF'$ we also have $\Na{\fF} \le \Na{\fF'}$.

\begin{lemm}\label{l:small support}
For every family $\fF$ of subsets of $\CC$ we have
$$
\sum_{W \in \fF} \diam(W)^\alpha \le 4^\alpha \cdot \diam(\supp(\fF))^\alpha \ \Na{\fF}.
$$
\end{lemm}
\begin{proof}
Recall that we identify the Riemann sphere $\CC$ with $\C \cup \{ \infty \}$.

When \, $\diam(\supp(\fF)) \ge \tfrac{1}{4}$ \, there is nothing to prove, so assume that \, $\diam(\supp(\fF)) < \tfrac{1}{4}$.
After an isometric change of coordinates assume $0 \in \supp(\fF)$.
The hypothesis \, $\diam(\supp(\fF)) < \tfrac{1}{4}$ implies that \, $\supp(\fF)$ \, is contained in $\{ z \in \C \mid |z| < 1 \}$.
Moreover, letting $\Delta \= 2\diam(\supp(\fF))$, we have that the set \, $\supp(\fF)$ \, is contained in the ball $B = \{ z \in \C \mid |z| \le \Delta \}$, because the density of the spherical metric with respect to the Euclidean metric on $\C$ is at least $\tfrac{1}{2}$ on $\{ z \in \C \mid |z| \le 1 \}$.

On the other hand, note that the distortion of the M\"obius map $\varphi(z) \= \Delta^{-1} z$ on $B$ is bounded by~$2$.
So we have
$$
\sum_{W \in \fF} \diam(W)^\alpha
\le
2^\alpha \Delta^{\alpha} \sum_{W \in \fF} \diam(\varphi(W))^\alpha
\le
4^\alpha \diam(\supp(\fF))^\alpha \Na{\fF}.
$$
\end{proof}
\subsection{Univalent pull-backs.}\label{ss:univalent pull-back}
The following assertion is an easy consequence of Koebe Distortion Theorem.
For every $\modul > 0$ there is a constant $C(\modul) > 0$ such that the following property holds.
Let $U$ and $\hU$ be simply-connected subsets of $\CC$ such that $\ov{U} \subset \hU$ and such that $\hU \setminus \ov{U}$ is an annulus of modulus at least $\modul$.
Then for every univalent holomorphic map $h : \hU \to \CC$ there exists a M\"obius map $\varphi$ such that the distortion of $h \circ \varphi$ on $U$ is bounded by $C(\modul)$.

The following lemma is a direct consequence of the property above.
\begin{lemm}\label{l:univalent pull-back}
For each $\modul > 0$ there is a constant $C_0(\modul) > 0$ such that the following properties hold.
Let $\hU$, $\hV$ be simply-connected subsets of $\CC$ and let $h : \hU \to \hV$ be a holomorphic map.
Let $\fF$ be a family of subsets of $\CC$ such that \ $\supp(\fF)$ \ is contained in a simply-connected subset $V$ of $\hV$, in such a way that $\hV \setminus \ov{V}$ is an annulus of modulus at least $\modul$.
Then,
$$
\Na{\{ h^{-1}(W) \mid W \in \fF \}} \le C_0(\modul) \Na{\fF}.
$$
\end{lemm}
\subsection{Modulus and diameter.}
Let $A$ be an open subset of $\CC$ homeomorphic to an annulus.
When the complement of $A$ in $\CC$ contains at least~$3$ points, there is a unique $R \in (1, + \infty]$ such that $A$ is conformally equivalent to $\{ z \in \C \mid 1 < |z| < R \} \in (0, + \infty]$.
In this case we put $\modulus(A) = \ln R$.
When the complement of $A$ in $\CC$ consists of two points, we put $\ln A = + \infty$.

A {\it round annulus} is by definition the complement in $\CC$ of~$2$ disjoint closed balls or the complement of a closed ball in an open ball containing it, that is different from $\CC$.
\begin{lemm}\label{l:small product}
There is a universal constant $\modul_0 > 0$ such that for every subset $A$ of $\CC$ that is homeomorphic to an annulus, the connected components $B$ and $B'$ of $\CC \setminus A$ satisfy
$$
\diam(B) \cdot \diam(B') \le \exp(-(\modulus(A) - \modul_0)).
$$
\end{lemm}
\begin{proof}
Is easy to check by direct computation that there is a constant $\modul_1 > 0$ so that property above holds for every round annulus, with $\modul_0$ replaced by $\modul_1$.
On the other hand, there is a universal constant $\modul_2 > 0$ such that every annulus in $\CC$ of modulus $\modul \in (\modul_1, + \infty)$ contains an essential round annulus of modulus $\modul - \modul_1$.
So the assertion of the lemma holds with constant $\modul_0 \= \modul_1 + \modul_2$.
\end{proof}

\subsection{Unicritical pull-backs.}\label{ss:unicritical pull-backs}
Given simply-connected subsets $U$ and $V$ of $\CC$, we say that a holomorphic and proper map $h : U \to V$ is {\it unicritical} if is has a unique critical point.
Note that the degree of a unicritical map as a ramified covering is equal to the local degree at its unique critical point.
Given an unicritical map $h : U \to V$ and a family $\fF$ of connected subsets of $\CC$ such that $\supp(\fF) \subset V$, put
$$
h^{-1}(\fF)
\=
\{ \text{connected component of $h^{-1}(W)$, for some $W \in \fF$} \}.
$$

We say that a family $\fF$ of simply-connected subsets of $\CC$ is $\modul$-{\it shielded} by a family $\hfF$, if for every $W \in \fF$ there is an element $\hW$ of $\hfF$ containing $W$ and such that $\hW \setminus \ov{W}$ is an annulus of modulus at least $\modul$.
\begin{lemm}\label{l:unicritical pull-back}
Given $\modul > 0$ and $d \ge 2$ there is a constant $C_1(\modul, d) > 0$ such that the following properties hold.
Let \, $h : \hU \to \hV$ be a unicritical map with critical point $c$ and let $\fF$ be a family of simply-connected subsets of $\CC$ that is $\modul$-shielded by a family $\hfF$ such that \, $h(c) \not \in \supp(\hfF)$ and such that \, $\supp(\hfF)$ \, is contained in a simply-connected subset $V$ of $\hV$, so that $\hV \setminus \ov{V}$ is a annulus of modulus at least $\modul$.
Then
$$
\Na{h^{-1}(\fF)} \le C_1(\modul, d) \cdot \Na{\fF}.
$$
\end{lemm}
\begin{proof}
Let $\gamma$ be the Jordan curve that divides the annulus $\hV \setminus \ov{V}$ into two annuli of modulus equal to a half of the modulus of $\hV \setminus \ov{V}$.
Denote by $\tV$ the open disk in $\hV$ bounded by $\gamma$.
Then there are two cases.

\partn{Case 1}
$h(c) \not \in \tV$.
Then the preimage of $\tV$ by $h$ has $d$ connected components and on each of these $h$ is univalent.
As by definition of $\tV$ the set $\tV \setminus \ov{V}$ is an annulus of modulus at least $\tfrac{1}{2}\modul$, in this case the assertion follows from Lemma~\ref{l:univalent pull-back}, with constant $C_1(\modul, d) \= d \cdot C_0(\tfrac{1}{2}\modul)$.

\partn{Case 2}
$h(c) \in \tV$.
In this case the set $\tU \= \varphi^{-1}(\tV)$ is connected and simply-connected.
Let $\varphi : \D \to \tV$ be a conformal uniformization such that $\varphi(0) = h(c)$.
As the modulus of the annulus $\tV \setminus \ov{V}$ is at least $\tfrac{1}{2}\modul$, there is $r_0 \in (0, 1)$ only depending on $\modul$ such that $V$ is contained in $\varphi(\{ |z| \le r_0 \})$.
By hypothesis for each $W \in \fF$ there is $\hW \in \hfF$ contained in $V \setminus \{ h(c) \}$, such that $\hW \setminus \ov{W}$ is an annulus of modulus at least $\modul$.
It follows that the diameter of $W$ with respect to the hyperbolic metric of $V \setminus \{ h(c) \}$ is bounded from above in terms of $\modul$ only.
Therefore there is $\rho \in (0, 1)$ such that for every $r \in (0, \rho r_0)$, every element $W$ of $\fF$ intersecting $\varphi(\{ |z| = r \})$ must be contained in the annulus $\varphi(\{ \rho r \le |z| \le \rho^{-1}r \})$.
For $j \ge 1$ put
$$
\tV_j \= \varphi (\{ r_0 \rho^{j + 1} \le |z| \le r_0 \rho^{j - 1} \})
\ \text{ and } \
\fF_j \= \{ W \in \fF \mid W \subset \tV_j \}.
$$
By definition of $\rho$ we have $\fF = \cup_{j \ge 1} \fF_j$.
Moreover put $\tU_j \= h^{-1}(\tV_j)$.
By~\S\ref{ss:univalent pull-back} there is a constant $C' > 0$ only depending on $r_0$ and on the degree $d$ of $h$, such that for every $j \ge 1$ we have
\begin{equation}\label{e:annulus Koebe}
\Na{h^{-1}(\fF_j)} \le C' \Na{ \fF_j} \le C' \Na{\fF}.
\end{equation}

Let $j_0 \ge 1$ be the least integer such that for every $j \ge j_0$ we have $\diam(\tU_j) < \tfrac{1}{2}$.
Note that for every $j > j_0$ the set
$$
h^{-1} \circ \varphi (\{ z \in \CC \mid r_0 \rho^{j - 1} < |z| < r_0 \rho^{j_0 - 1} \})
$$
is an annulus of modulus $(j - j_0) \ln \rho^{-1} / d$.
So there are constants $C'' > 0$ and $\eta \in (0, 1)$ only depending on $\rho$ and $d$, such that for every $j \ge j_0$ we have $\diam(\tU_j) \le C'' \eta^{(j - j_0)}$.
Therefore, letting $C''' \= 4^\alpha C'C'' (1 - \eta)^{-1}$, Lemma~\ref{l:small support} implies that,
$$
\sum_{W \in \cup_{j \ge j_0} h^{-1}(\fF_j)} \diam(W)^\alpha \le C''' \Na{\fF}.
$$

Let $\modul_0$ be the constant given by Lemma~\ref{l:small product} and let $N$ be the smallest positive integer such that $\exp(-(N d^{-1} \ln \rho^{-1} - \modul_0)) < \tfrac{1}{4}$.
Note that $N$ depends on $\modul$ only.
By definition of $j_0$ we have $\diam(\tU_{j_0 - 1}) \ge \tfrac{1}{2}$.
As for every $j \le j_0 - (N + 3)$ the sets $\tU_j$ and $\tU_{j_0 - 1}$ are separated by an annulus of modulus $(j_0 - j - 3) d^{-1} \ln \rho^{-1} \ge N d^{-1} \ln \rho^{-1}$, Lemma~\ref{l:small product} implies that $\diam(\tU_j) < \tfrac{1}{2}$ and as before we have,
$$
\sum_{W \in \cup_{j = 1, \ldots, j_0 - (N + 3)} h^{-1}(\fF_j)} \diam(W)^\alpha \le C''' \Na{\fF}.
$$
Therefore,
\begin{multline*}
\sum_{W \in h^{-1}(\fF)} \diam(W)^\alpha
\le 
\sum_{j = j_0 - (N + 2), \ldots, j_0 - 1} \Na{h^{-1}(\fF_j)} + 2 C''' \Na{\fF}
\\ \le
(C'(N + 2) + 2C''') \Na{\fF}.
\end{multline*}
As our hypothesis are coordinate free, this estimate proves the assertion of the lemma with constant $C_1(\modul, d) \= C' (N + 2) + 2 C'''$.
\end{proof}

\section{Nice sets and the density $\Na{\cdot}$.}\label{s:nice sets and the density}
Fix a complex rational map $f$ of degree at least~$2$.
Recall that the {\it hyperbolic Hausdorff dimension} of $f$ is by definition
$$
\HDhyp(f) \= \sup_{X} \HD(X),
$$
where the supremum is taken over all forward invariant subsets $X$ of $\CC$ on which $f$ is uniformly expanding.

The purpose of this section is to prove the following proposition.
\begin{prop}\label{p:finite pressure}
Let $f$ be a rational map such that for every neighborhood $V'$ of $\sCJ$ the map $f$ is uniformly expanding on the set
\begin{equation}\label{e:points avoiding V'}
\{ z \in J(f) \mid \text{ for every $n \ge 0$ we have $f^n(z) \not \in V'$} \}.
\end{equation}
Then for every nice couple $(\hV, V)$ for $f$ there exists $\alpha \in (0, \HDhyp (f))$ such that the collection $\fD_V$ of connected components of $\CC \setminus K(V)$ satisfies $\Na{\fD_V} < + \infty$.
\end{prop}
See~\S\ref{ss:nice sets} for the definition of $K(V)$.

After the preliminary lemmas~\ref{l:hd}, \ref{l:Whitney distribution} and~\ref{l:sparse}, the proof of this proposition is given at the end of this subsection.

For a subset $X$ of $\CC$ denote by $\BD(X)$ the upper box dimension of $X$.
\begin{lemm}\label{l:hd}
  Let $f$ be a rational map as in the statement of Proposition~\ref{p:finite pressure}.
If $K'$ is a compact subset of $\CC$ such that $f(K') \subset K'$ and such that $f$ is uniformly expanding on $K'$, then
$$
\HD(K') \le \BD(K') < \HDhyp(J(f)).
$$
\end{lemm}
\begin{proof}
We will show first that we can reduce to the case when $K'$ is compact and forward invariant ($f(K') = K'$).
Let $K$ be the closure of $\cap_{n = 0}^{\infty} f^n(K')$.
Clearly $K$ is compact and forward invariant.
Thus there is a compact neighborhood $U$ of $K$ such that $U$ is contained in the interior of $f(U)$.
Let $\tK$ be the maximal invariant set of~$f$ contained in $U$.
This set is compact, forward invariant and, reducing $U$ if necessary, $f$ is uniformly expanding on $\tK$.
By definition of $K$ it follows that there is a sufficiently large integer~$N$ such that $f^N(K')$ is contained in~$U$, and hence in $\tK$.
It follows that $\BD(K') \le \BD(\tK)$.
Thus, replacing $K'$ by $\tK$ if necessary, we can reduce to the case in which $K'$ is compact and invariant.

The first inequality is a general fact.
As~$f$ is uniformly expanding on~$K'$, it follows that there is
a uniformly expanding set $K''$ containing $K'$, which has a Markov partition,
see [PU].
(The idea, taken from Bowen's construction, is to pick a $\delta$-net in $K'$ and shadow by forward $f$-trajectories all $\epsilon$-trajectories of points in the net. The set $K''$ is defined as the union of these trajectories and it contains $K'$, provided $\delta \ll \epsilon$.)
Hence $\BD(K'')=\HD(K'')$.
To get this one uses an equilibrium state
for the potential $-\HD(K'')\ln |f'|$ on an invariant topologically
transitive part of the related topological Markov shift.
Next, for a sufficiently small neighborhood $V'$ of $\sC$ one can construct in a similar way an expanding set $K'''$ containing the set~\eqref{e:points avoiding V'}, and that has a Markov partition that is a substantial extension of the previous one.
Hence we have
$$
\BD(K') \le \BD(K'') = \HD(K'') < \HD(K''') \le \HDhyp(J(f)).
$$
\end{proof}
\begin{lemm}\label{l:Whitney distribution}
Let $f$, $(\hV, V)$ and $\fD_V$ be as in Proposition~\ref{p:finite pressure}.
Then there are constants $\alpha \in (0, \HDhyp(f))$ and $C_0 > 0$ such that
\begin{equation}\label{e:total sum}
\sum_{W \in \fD_V} \diam(W)^\alpha < + \infty,
\end{equation}
and such that for every ball $B$ of $\CC$ we have,
\begin{equation}\label{e:sum for balls}
\sum_{W \in \fD_V, \, W \subset B} \diam(W)^\alpha \le C_0 \diam(B)^\alpha.
\end{equation}
\end{lemm}
\begin{proof}
As by hypothesis $f$ is uniformly expanding on $K \= K(V) \cap J(f)$, it follows that for every ball $B$ of $\CC$ intersecting $K$ there in an integer $n \ge 0$ such that $f^n(B)$ has definite size and such that the distortion of $f^n$ on $B$ is bounded independently of $B$.
So the existence of the constant $C_0 > 0$ for which~\eqref{e:sum for balls} holds for every ball $B$ follows from~\eqref{e:total sum}.
In what follows we will prove that~\eqref{e:total sum} holds for some appropriated choice of $\alpha$.

\partn{1}
Let $V'$ be a neighborhood of $\sCJ$ that is contained in $V$ and such that for every $c \in \sCJ$ the set $K'$ defined by~\eqref{e:points avoiding V'} intersects $V^c$. 
It follows that every element of $\fD_V$ intersects $K'$.
Given $W \in \fD_V$ choose a point $z_W$ in $W \cap K'$ and put $B_W \= B(W, \diam(W))$, so that $W \subset B_W$.

By Lemma~\ref{l:hd} we have $\BD(K') < \HDhyp(f)$.

\partn{2}
Fix $\alpha \in (\BD(K'), \HDhyp(f))$.
For $r > 0$ let $N(r)$ be the least number of balls of $\CC$ of radius $r$ centered at some point of $K'$, that are necessary to cover~$K'$.
It follows that for every $r > 0$ the cardinality of a collection of pairwise disjoint balls of radius at least $r$ and centered at points in $K'$ is at most $N(r)$.
On the other hand, as
$$
\alpha > \BD(K') = \limsup_{r \to 0} \frac{\ln N(r)}{\ln(r^{-1})},
$$
it follows that,
\begin{equation}\label{e:total box sum}
\sum_{n \ge 0} N(2^{-n}) 2^{-n\alpha} < + \infty.
\end{equation}

\partn{3}
Since for every $W \in \fD_V$ the map $f^{m_W}$ extends univalently to $\hW$, it follows that there is a constant $A_1 > 0$ such that for every $W \in \fD_V$ we have $\Area(W) \ge A_1 \diam(W)^2$.
Therefore there is an integer $A_2$ such that for every $W \in \fD_V$ there are at most $A_2$ elements $W' \in \fD_V$ such that,
$$
\tfrac{1}{2} \diam(W) \le \diam(W') \le 2 \diam(W),
$$
and such that $B_{W'}$ intersects $B_W$.

So, if for $r > 0$ we denote by $N'(r)$ the number of elements of $\fD_V$ whose diameter belongs to $(r, 2r)$, then we can find a collection $\fF$ of such sets whose cardinality is at least $N'(r) / (A_2  + 1)$ and such that the balls $B_W$, for $W \in \fF$, are pairwise disjoint.
As the centers of the balls $B_W$ belong to $K'$, it follows that
$$
N'(r) \le (A_2 + 1) N(r).
$$
So by~\eqref{e:total box sum} we have,
$$
\sum_{W \in \fD_V} \diam(W)^\alpha
\le
\sum_{n \ge 0} N'(2^{-n}) 2^{-n\alpha} < + \infty.
$$
\end{proof}

\begin{lemm}\label{l:sparse}
Let $f$ be a complex rational map and let $(\hV, V)$ be a nice couple for $f$.
Then there is a constant $C_1 > 0$ such that for every ball $B$ of $\CC$ there is at most one connected component $W$ of $\CC \setminus K(V)$ intersecting $B$ and such that
$$
\diam(W) \ge C_1 \diam(B).
$$
\end{lemm}
\begin{proof}
Let $(\hV, V)$ be a nice couple for $f$.
It follows from~\eqref{e:Markov property for pull-backs} that if $W$ and $W'$ are distinct elements of $\fD_V$, then we have either $W \cap \hW' = \emptyset$ or $W \cap \widehat{W'} = \emptyset$.
So to prove the lemma is enough to show that there is a constant $C_1 > 0$ such that if $W$ is an element of $\fD_V$ and if $B$ is a ball intersecting $W$ such that $\diam(W) \ge C_1 \diam(B)$, then $B \subset \hW$.

Note that for each $W \in \fD_V$ the set $\hW$ is disjoint from $K(\hV)$ and $\hW \setminus \ov{W}$ is an annulus of modulus at least
$$
\modul \= \min_{c \in \sCJ} \modulus(\hV^c \setminus \ov{V^c}).
$$
As the set $K(\hV)$ contains at least~$2$ points (because it contains $\partial V$), it follows that there is $\varepsilon > 0$ such that for every $W \in \fD_V$ the set
$$
W(\varepsilon) = \{ z \in \CC \mid \dist(z, W) \le \varepsilon \diam(W) \}
$$
is contained in $\hW$.
Clearly the constant $C_1 \= \varepsilon^{-1}$ has the desired property.
\end{proof}
\begin{proof}[Proof of Proposition~\ref{p:finite pressure}]
Recall that we identify $\CC$ with $\C \cup \{ \infty \}$ and that the spherical metric is normalized in such a way that its density with respect to the Euclidean metric on $\C$ is given by $z \mapsto (1 + |z|^2)^{-1}$.
Moreover, balls, distances, diameters and derivatives are all taken with respect to the spherical metric, see the Preliminaries.

\medskip

Let $\alpha$, $C_0$ and $C_1$ be the constants given by lemmas~\ref{l:Whitney distribution} and~\ref{l:sparse}.
Choose $\rho \in (0, 1)$ sufficiently small so that for every $r \in (0, 1]$ we have
$$
\dist(\{ |z| = r \}, \{ |z| \le \rho r \})
\ge
C_1 \diam( \{ |z| \le \rho r \}).
$$
Let $\varphi$ be a given M\"obius map.
After isometric change of coordinates in the domain and in the target, we assume that $\varphi$ is of the form $\varphi(z) = \lambda z$, with $\lambda$ real and satisfying $\lambda \ge 1$.
Let $N \ge 1$ be the least integer satisfying $\rho^N < \lambda^{- 1/2}$ and put
$$
A_0 \= \{ |z| > \rho \lambda^{-1/2} \} \cup \{ \infty \},
$$
$$
A_n \= \{ \rho^{n + 1} \lambda^{-1/2} < |z| < \rho^{n - 1} \lambda^{-1/2} \},
$$
for $n = 1, \ldots, N - 1$, and put
$$
A_N \= \{ |z| < \rho^{N - 1} \lambda^{-1/2} \}.
$$
For $n = 0, \ldots, N$ we denote by $\fF_n$ the sub-collection of $\fD_V$ of sets contained in $A_n$ and we denote by $\fF$ the sub-collection of $\fD_V$ of sets not contained in any of the $A_n$.
By definition we have $\fF \cup \fF_0 \cup \ldots \cup \fF_N = \fD_V$.

In parts~$1$, $2$ and~$3$ below we estimate the sum $\sum \diam(\varphi(W))^\alpha$, where $W$ runs through $\fF$, $\fF_0$ and $\fF_1 \cup \ldots \cup \fF_N$, respectively.
These estimates are independent of $\lambda$, so the lemma follows form them.

\partn{1}
By definition each element of $\fF$ intersects at least two of the sets $\{ |z| = \rho^n \lambda^{-1/2} \}$, for $n = 0, \ldots, N$.
For $W \in \fF$ denote by $n(W)$ the largest integer $n = 0, \ldots, N$ such that $W$ intersects $\{ |z| = \rho^n \lambda^{-1/2} \}$.
So for each $W \in \fF$ we have
$$
\varphi(W)
\subset
\varphi(\{ |z| \ge \rho^{n(W) + 1} \lambda^{-1/2} \} \cup \{ \infty \} )
\subset \{ |z| \ge \rho^{n(W) + 2 - N} \} \cup \{ \infty \},
$$
and $\diam(\varphi(W)) \le (2\rho^{-2}) \rho^{n(W) - N}$.
By Lemma~\ref{l:sparse} and by the choice of $\rho$, for distinct $W, W' \in \fF$ the integers $n(W)$ and $n(W')$ are distinct.
Therefore we have
\begin{multline*}
\sum_{W \in \fF} \diam(\varphi(W))^\alpha
\le
(2\rho^{-2})^\alpha (1 + \rho^\alpha + \ldots + \rho^{N \alpha})
\\ \le
(2\rho^{-2})^\alpha/(1 - \rho^\alpha).
\end{multline*}

\partn{2}
As $\sup \{ |\varphi'(z)| \mid z \in A_0 \} \le \rho^{-2}$, we have
$$
\sum_{W \in \fF_0} \diam(\varphi(W))^\alpha
\le
\rho^{- 2\alpha} \sum_{W \in \fF_0} \diam(W)^\alpha
\le
\rho^{- 2\alpha} \sum_{W \in \fD_V} \diam(W)^\alpha.
$$
By Lemma~\ref{l:Whitney distribution} this last quantity is finite.

\partn{3}
Note that there is a constant $C_2 > 0$ only depending on $\rho$, that bounds the distortion of $\varphi$ on each of the sets $A_n$, for $n = 1, \ldots, N$.
On the other hand, note that for each $n = 1, \ldots, N$ and $z \in \C$ such that $|z| = \rho^n \lambda^{-1/2}$, we have $z \in A_n$ and
$$
|\varphi'(z)|
=
\lambda \frac{1 + \rho^{2n} \lambda^{-1}}{1 + \rho^{2n} \lambda}
\le 2 \rho^{-2n}.
$$
Therefore we have
$$
\sum_{W \in \fF_n} \diam(\varphi(W))^\alpha
\le
(2C_2)^\alpha \rho^{-2n\alpha} \sum_{W \in \fF_n} \diam(W)^\alpha.
$$
Since $\supp(\fF_n) \subset A_n \subset \{ |z| \le \rho^{n - 1} \lambda^{-1/2} \} \subset \{ |z| \le \rho^{n + N - 2} \}$, so
$$
\sum_{W \in \fF_n} \diam(W)^\alpha
\le
C_0 \diam(A_n)^\alpha
\le
C_0 (2\rho^{-2})^\alpha \rho^{(n + N) \alpha}.
$$
So, letting $C_3 \= (2C_2)^\alpha C_0 (2 \rho^{-2})^\alpha$ we have
\begin{multline*}
\sum_{W \in \fF_n} \diam(\varphi(W))^\alpha
\le C_3 \rho^{(N - n)\alpha}
\ \text{ and} \\
\sum_{W \in \fF_1 \cup \ldots \fF_N} \diam(\varphi(W))^\alpha \le C_3/(1 - \rho^\alpha).
\end{multline*}
\end{proof}

\section{Key Lemma.}\label{s:key lemma}
The purpose of this section is to prove the following lemma.
\begin{generic}[Key Lemma]
Let $f$ be a rational map satisfying \ExpShrink.
Then for every $\modul > 0$ there exists $r > 0$ such that if $(\hV, V)$ is a nice couple satisfying,
\begin{equation}\label{e:strong regularity hypothesis}
\min_{c \in \sCJ} \modulus (\hV^c \setminus \ov{V^c}) \ge \modul
\ \text{ and } \
\max_{c \in \sCJ} \diam(\hV^c) \le r,
\end{equation}
then the canonical induced map $F : D \to V$ associated to $(\hV, V)$ satisfies the following properties.
\begin{enumerate}
\item[1.]
We have $\HD((J(f) \setminus J(F)) \cap V) < \HD(J(f))$.
In particular, for every $c \in \sC$ we have $\HD(J(F) \cap V^c) = \HD(J(f))$.
\item[2.]
There exists $\alpha \in (0, \HD(J(f)))$ such that
\begin{equation*}
\sum_{W \text{ c.c. of } D} \diam(W)^\alpha < + \infty,
\end{equation*}
where the sum is over all connected components $W$ of $D$.
\end{enumerate}
\end{generic}

After some preliminarily considerations in~\S\ref{ss:bad pull-backs}, we prove parts~$1$ and~$2$ of the Key Lemma in~\S\ref{ss:Key Lemma 1} and in~\S\ref{ss:Key Lemma 2}, respectively.
\subsection{Bad pull-backs.}\label{ss:bad pull-backs}
Let $f$ be a rational map and let $\hV$ be a nice set for $f$.
For an integer $n \ge 1$ we will say that a connected component $\tW$ of $f^{-n}(\hV)$ is a {\it bad pull-back of $\hV$ of order} $n$, if $f^n$ is not univalent on $\tW$ and if for every $m = 1, \ldots, n - 1$ such that $f^m(\tW) \subset V$, the map $f^m$ is not univalent on the connected component of $f^{-m}(\hV)$ containing $\tW$.
Note that every bad pull-back of $\hV$ contains a critical point of $f$ in~$\sCJ$.

The proof of the following lemma is similar to that of Lemma~$A.2$ of~\cite{PRS}.
\begin{lemm}\label{l:bad pull-backs}
Let $f$ be a rational map and let $\hV$ be a nice set for $f$.
Let $L \ge 1$ be the least integer such that for some $c \in \sCJ$ we have $f^L(c) \in \hV$.
Then for each positive integer $n$, there are at most $(2 L \#\sCJ)^{2n/L}$ bad pull-backs of $\hV$ of order~$n$.
\end{lemm}
\begin{proof}

\partn{1}
For a bad pull-back $\tW$ of $\hV$ of order $n$, let $\ell(\tW)$ be the largest integer $\ell$ in $\{0, \ldots, n - 1 \}$ such that $f^\ell(\tW)$ intersects $\sCJ$.
The integer $\ell(\tW)$ might be equal to~$0$.
As $\hV$ is a nice set we have $f^{\ell(\tW)}(\tW) \subset \hV$ and, as $W$ is a bad pull-back, the connected component of $f^{-\ell(\tW)}(\hV)$ containing $\tW$ is a bad pull-back of $\hV$ of order $\ell(\tW)$.

\partn{2}
Fix an integer $n \ge 1$.
For a given bad pull-back $W$ of order $n$ define a strictly decreasing sequence of integers $(\ell_0, \ldots, \ell_k)$, by induction as follows.
Define $\ell_0 \= n$ and suppose that for some $j \ge 0$ the integer $\ell_j$ is already defined and that we have $f^{\ell_j}(\tW) \subset \hV$.
If $\ell_j = 0$, then define $k \= j$ and stop.
If $\ell_j > 0$, then by induction hypothesis $f^{\ell_j}(\tW)$ is contained in $\hV$ and therefore the connected component $\tW'$ of $f^{-\ell_j}(\hV)$ containing $\tW$ is a bad pull-back of $\hV$ of order $\ell_j$.
Then define $\ell_{j + 1} \= \ell(\tW')$.
As remarked in part~$1$, in this case we have $f^{\ell_{j + 1}}(\tW) \subset f^{\ell_{j + 1}}(\tW') \subset \hV$, so the induction hypothesis is satisfied.

\partn{3}
To each bad pull-back of $\hV$ of order $n$ we associate a strictly decreasing sequence $(\ell_0, \ldots, \ell_k)$, as in part~$2$, so that $\ell_0 = n$ and $\ell_k = 0$.
Note that for each $j = 1, \ldots, k$, the pull-back of $\hV$ by $f^{\ell_{j - 1} - \ell_{j}}$ containing $f^{\ell_{j}}(\tW)$ contains a critical point in $\sCJ$.
As for each $c \in \sCJ$ and each integer $m$ there are at most $\# \sCJ$ connected components of $f^{-m}(\hV^c)$ containing an element of $\sCJ$, it follows that there are at most $(\# \sCJ)^{k + 1}$ bad pull-backs of order $n$ with the same associated sequence.

On the other hand, by definition of $L$ it follows that for every $j = 1, \ldots, k$ we have $\ell_{j - 1} - \ell_{j} \ge L$.
So, $k \le n/L$ and for each integer $m = 1, \ldots, n$ there is at most one integer $r \in \{0, 1, \ldots, L - 1 \}$ such that $m + r$ is one of the $\ell_j$.
It follows that there are at most $(L + 1)^{2n/L}$ such decreasing sequences.

We conclude that the number of bad pull-backs of $\hV$ of order~$n$ is at most,
$$
(\# \sCJ)^{k + 1} (L + 1)^{2n/L} \le (2 L \# \sCJ)^{2n/L}.
$$
\end{proof}

\subsection{Proof of part~$1$ of the Key Lemma.}\label{ss:Key Lemma 1}
Let $f$ be a rational map satisfying condition \ExpShrink{} with constants $\lambda_{\Exp} > 1$ and $r_0 > 0$.
We will show that part~$1$ of the Key Lemma holds for every nice couple $(\hV, V)$ for which $\max_{c \in \sCJ} \diam(V^c)$ is sufficiently small.
In fact, we will prove that for such a nice couple $(\hV, V)$, we have
\begin{equation}\label{e:hd estimate}
\HD((V \cap J(f)) \setminus D) < \HD(J(f)).
\end{equation}
As every inverse branch of $F$ is Lipschitz, this implies that $\HD((J(f) \setminus J(F)) \cap V) < \HD(J(f))$, as desired.

Given a nice couple $(\hV, V)$, denote by $R_V$ the first return map to $V$ and denote by $J(R_V)$ the subset of $V$ of those points for which $R_V^n$ is defined for every integer $n \ge 1$.
Equivalently, $J(R_V)$ is the set of those points that return infinitely often to~$V$ under forward iteration of $f$.
Note that $J(F) \subset J(R_V)$ and that
$$
\HD((V \cap J(f)) \setminus J(R_V)) \le \HD(K(V) \cap J(f)).
$$
As $f$ is uniformly expanding on $K(V) \cap J(f)$ (this follows easily from condition \ExpShrink), Lemma~\ref{l:hd} implies that,
$$
\HD((V \cap J(f)) \setminus J(R_V)) < \HD(J(f)).
$$
So, the following lemma implies that~\eqref{e:hd estimate} holds for every nice couple $(\hV, V)$ for which $\max_{c \in \sCJ} \diam(V^c)$ is sufficiently small, by choosing $\varepsilon \in (0, \HD(J(f)))$.

For a given $r > 0$ let $L(r) \ge 1$ be the smallest integer such that for some $c \in \sCJ$ the point $f^{L(r)}(c)$ is at distance at most $r$ from $\sCJ$.
As our standing convention is that no critical point in $\sCJ$ is mapped to a critical point under forward iteration, we have that $L(r) \rightarrow \infty$ as $r \rightarrow 0$.
\begin{lemm}\label{l:HD of bad}
Given $\e > 0$ choose $r \in (0, r_0)$ sufficiently small so that
\begin{equation}\label{e:Hausdorff}
(2L(r)\#\sCJ)^{2/L(r)}\lambda_{\Exp}^{-\e} < 1.
\end{equation}
Then for every nice couple $(\hV, V)$ such that $\max_{c \in \sCJ} \diam(\hV^c) \le r$, the Hausdorff dimension of $J(R_V) \setminus D$ is at most $\varepsilon$.
\end{lemm}
\begin{proof}
For a point $z$ in $J(R_V) \setminus D$ there are arbitrarily large integers $n \ge 1$ such that $f^n(z) \in V$.
Moreover, for every such $n$ the pull-back $\tW$ of $\hV$ to $z$ by $f^n$ is not univalent.
It follows that $\tW$ is a bad pull-back of $\hV$.
So, if for $n \ge 1$ we denote by $\fD_n$ the collection of all bad pull-backs of $\hV$ of order~$n$, then for every $n_0 \ge 1$ we have
$$
J(R_V) \setminus D \subset \cup_{n \ge n_0} \cup_{\tW \in \fD_n} \tW.
$$
As,
$$
\sum_{\tW \in \fD_n} \diam(\tW)^\e
\le
(2L(r)\# \sCJ)^{2n/L(r)} \lambda_{\Exp}^{- n\varepsilon},
$$
it follows from \eqref{e:Hausdorff} that this sum is exponentially small with~$n$.
This implies the assertion of the lemma.
\end{proof}
\subsection{Proof of part~$2$ of the Key Lemma.}\label{ss:Key Lemma 2}
Let $f$, $\la_{\Exp}$, $r_0 > 0$ and $L(r)$ be as in the previous subsection.
For $c \in \sCJ$ let $d(c) > 1$ be the local degree of $f$ at $c$.
Let $\modul > 0$ be given and for an integer $d \ge 2$ let $C_1(\modul, d) > 0$ be the constant given by Lemma~\ref{l:unicritical pull-back} and put
$$
C_1(\modul) \= \max_{c \in \sCJ} C_1(\modul, d(c)).
$$
Let $r \in (0, r_0)$ be sufficiently small so that
\begin{equation}\label{e:strong regularity}
(2L(r)(\#\sCJ))^{2/L(r)} \lambda_{\Exp}^{- \HD(J(f))} C_1(\modul)^{1/L(r)} 
<
1.
\end{equation}
We will prove that part~$2$ of the Key Lemma holds for this choice of $r$.
So let $(\hV, V)$ be a nice couple satisfying,
$$
\min_{c \in \sCJ} \modulus (\hV^c \setminus \ov{V^c}) \ge \modul
\ \text{ and } \
\max_{c \in \sCJ} \diam(\hV^c) \le r.
$$

Denote by $\fD_V$ the collection of all connected components of $\CC \setminus K(V)$ and let $\alpha \in (0, \HDhyp(J(f)))$ be given by Proposition~\ref{p:finite pressure}, so that $\Na{\fD_V} < + \infty$ (it is easy to see that the hypothesis of this proposition are satisfied for maps satisfying property \ExpShrink).
Taking $\alpha$ closer to $\HD(J(f)) = \HDhyp(J(f))$ if necessary, we assume that~\eqref{e:strong regularity} holds with $\HD(J(f))$ replaced by~$\alpha$.

Recall that we denote by $D$ the subset of $V$ of those points having a good time, see~\S\ref{ss:canonical induced}.
Denote by $\fD$ the collection of the connected components of~$D$.
As $V$ is a nice set, it is easy to see that every $W \in \fD$ is a pull-back of $V$ and that for every $z \in W$ we have $m(z) = m_W$.
\begin{lemm}\label{l:good at first}
Let $\fD_0$ be the sub-collection of $\fD$ of those $W$ such that $m_W$ is the first return time of~$W$ to $V$.
Then $\Na{\fD_0} < + \infty$.
\end{lemm}
\begin{proof}
Let $\fD_V'$ be the sub-collection of $\fD_V$ of sets $W'$ in $\fD_V$ such that $\widehat{W'}$ is disjoint from the critical values of $f$.
We have $\Na{\fD_V'} \le \Na{\fD_V} < + \infty$ and $\fD_0 \subset f^{-1}(\fD_V')$.
Then the assertion of the lemma follows from Lemma~\ref{l:unicritical pull-back} applied to $h \= f$.
\end{proof}
\begin{lemm}\label{l:good in bad}
\

\begin{enumerate}
\item[1.]
Every $W \in \fD \setminus \fD_0$ is contained in a bad pull-back $\tW$ of $\hV$ such that $m_{\tW} < m_W$ and that $f^{m_{\tW}}(W) \in \fD_0$.
\item[2.]
For every bad pull-back $\tW$ of $\hV$, the collection $\fD_{\tW}$ of all $W \in \fD$ contained in $\tW$ and such that $f^{m_{\tW}}(W) \in \fD_0$, satisfies,
$$
\Na{\fD_{\tW}} \le C_1(\modul)^{m_{\tW}/L(r)} \Na{\fD_0}.
$$
\end{enumerate}
\end{lemm}
\begin{proof}

\partn{1}
Let $W \in \fD \setminus \fD_0$ be given and let $n$ be the largest integer in $\{ 0, \ldots, m_W - 1\}$ such that $f^n(W)$ intersects $V$.
We have $n > 0$ because by assumption $W \not \in \fD_0$.
Then $W' \= f^n(W) \subset V$ and, as $f^{m_W}$ is univalent on $\hW$, we have that $f^{m_{W'}}$ is univalent on $\widehat{W'}$.
Moreover, by maximality of $n$ it follows that $f^{m_{W'}}$ coincides with the first return map to $V$ on $W'$, so that $W' \in \fD_0$.
To finish the proof is enough to show that the pull-back $\tW$ of $\hV$ by $f^n$ containing $W$ is bad.
If $\tW$ is not a bad pull-back of $\hV$, then there is $m \le m_{\tW} = n < m_W$ such that $f^m(W) \subset V$ and such that the pull-back of $\hV$ by $f^m$ containing $W$ is univalent.
But this contradicts the fact that $W$ is a connected component of $D$.

\partn{2}
We keep the notation of parts~$1$ and~$2$ of the proof of Lemma~\ref{l:bad pull-backs}.
Note that we have $\ell_0 = m_{\tW}$ and for each $j = 1, \ldots, k$ we have $\ell_{j - 1} - \ell_j \ge L(r)$, so that $k \le m_{\tW}/L(r)$. 
For each $j = 0, \ldots, k$ denote by $c_j$ the critical point such that $f^{\ell_j}(\tW) \subset \hV^{c_j}$ and let $\hU_j$ be the connected component of $f^{-(\ell_{j} - \ell_{j + 1})}(\hV^{c_{j}})$ containing $c_{j + 1}$.
It follows that $h_j : \hU_j \to \hV^{c_j}$ defined by $h_j \= f^{\ell_j - \ell_{j + 1}}|_{\hU_j}$ is a unicritical map with critical point $c_{j + 1}$ and degree $d(c_{j + 1})$.

Define families $\fF_0, \ldots, \fF_k$ of pull-backs of $V$ inductively as follows.
Put $\fF_0 \= \fD_0$ and for $j = 0, \ldots, k - 1$ suppose that the family $\fF_j$ is already defined and that $\supp(\fF_j) \subset f^{\ell_j}(\tW) \subset \hV^{c_j}$.
Let $\fF_j'$ be the sub-family of $\fF_j$ of those $W$ such that $\hW$ is disjoint from $h_j(c_{j + 1})$ and define
$$
\fF_{j + 1} \= h_j^{-1}(\fF_j').
$$
Clearly $\supp(\fF_{j + 1}) \subset \hU_j \subset V^{c_{j + 1}} \subset \hV^{c_{j + 1}}$, so the induction hypothesis is satisfied.

It is easy to see that $\fF_k = \fD_{\tW}$.
To prove the assertion of the lemma, note that the family $\fF_j'$ is $\modul$-shielded by the family $\hfF_j' = \{ \hW \mid W \in \fF_j' \}$ and that by definition of $\fF_{j + 1}$ we have $h(c_{j + 1}) \not \in \supp(\hfF_j')$. 
Moreover, by induction in $j$ we have that $\supp(\hfF_j') \subset V^{c_j}$.
So Lemma~\ref{l:unicritical pull-back} implies that
$$
\Na{\fF_{j + 1}}
\le
C_1(\modul, d(c_{j + 1})) \Na{\fF_j}
\le
C_1(\modul) \Na{\fF_j}.
$$
Therefore we have,
$$
\Na{\fD_{\tW}} = \Na{\fF_k} \le C_1(\modul)^k \Na{\fF_0} \le C_1(\modul)^{m_{\tW}/L(r)} \Na{\fD_0}.
$$
\end{proof}

\

To prove part~$2$ of the Key Lemma, note that by part~$1$ of Lemma~\ref{l:good in bad} we have
$$
\fD = \fD_0 \cup (\cup_{\tW} \fD_{\tW}),
$$
where the union is over all bad pull-backs $\tW$ of $\hV$.
Let $n \ge 1$ be an integer and let $\tW$ be a bad pull-back of order $n$.
Then we have $\diam(\tW) \le \la_{\Exp}^{-n}$.
By part~$2$ of Lemma~\ref{l:good in bad} and by Lemma~\ref{l:small support} we have
$$
\sum_{W \in \fD_{\tW}} \diam(W)^\alpha
\le
4^\alpha \lambda_{\Exp}^{-n \alpha} C_1(\modul)^{n / L(r)} \Na{ \fD_0 }.
$$
As the number of bad pull-backs of order $n$ is bounded by $(2L(r)(\# \sCJ))^{2n/L(r)}$, letting
$$
\eta = (2L(r)(\# \sCJ))^{2/L(r)}\lambda_{\Exp}^{- \alpha} C_1(\modul)^{1/L(r)},
$$
we have
$$
\sum_{\tW} \sum_{W \in \fD_{\tW}} \diam(W)^\alpha \le 4^\alpha \eta^n \Na{\fD_0},
$$
where the sum is over all bad pull-backs $\tW$ of $\hV$ of order~$n$.
Recall that we have chosen $r > 0$ sufficiently small and $\alpha$ sufficiently close to $\HD(J(f))$, so that~\eqref{e:strong regularity} is satisfied with $\HD(J(f))$ replaced by $\alpha$.
As $\alpha \in (0, \HD(J(f)))$ it follows that $\eta \in (0, 1)$ and that,
$$
\sum_{W \in \fD} \diam(W)^\alpha \le 4^\alpha(1 - \eta)^{-1} \Na{ \fD_0 } < + \infty.
$$

\section{Conformal and invariant measures.}\label{s:conformal and invariant}
In this section we prove the main results of this paper.
\subsection{Proof of theorem~\ref{t:conformal measures}.}\label{ss:proof of A}
Let $f$ be a rational map satisfying the \ExpShrink{} condition with constant $\lambda_{\Exp} > 1$, and let $\lambda \in (1, \la_{\Exp})$ and $\bfm > 0$ be given.
Let $r > 0$ be given by the Key Lemma for this choice of $\bfm$ and let $(\hV, V)$ be a nice couple for $f$ given by Theorem~\ref{t:induced map} for this choice of~$\lambda$, $\bfm$ and~$r$.
It follows that the nice couple $(\hV, V)$ satisfies the conclusions of the Key Lemma.
Reducing $r > 0$ if necessary, we assume that the canonical induced map associated to $(\hV, V)$ is topologically mixing (Lemma~\ref{l:mixingness}). 
Then Theorem~\ref{t:conformal measures} is a direct consequence of Theorem~\ref{t:conformal via inducing} in Appendix~\ref{a:conformal via inducing}, and of the fact that the equality~$\alpha(f) = \HD(J(f))$ holds for maps satisfying \ExpShrink, see~\cite{Prconical}.
\subsection{Proof of theorems~\ref{t:invariant measures} and~\ref{t:statistical properties}.}\label{ss:proof of B and C}
The uniqueness part of Theorem~\ref{t:invariant measures} follows from Proposition~\ref{p:conical conformal}, and from the fact that the unique conformal probability measure of minimal exponent of~$f$ is supported on the conical Julia set (Theorem~\ref{t:conformal measures}).
All the remaining statements will be obtained from some results of Young~\cite{You}, that we recall now.

Let $(\Delta_0, \cB_0, \measure_0)$ be a measurable space and let $T_0 : \Delta_0 \to \Delta_0$ be a measurable map for which there is a countable partition $\sP_0$ of $\Delta_0$, such that for each element $\Delta'$ of $\sP_0$ the map $T_0 : \Delta' \to \Delta_0$ is a bijection.
Moreover we assume that the partition $\sP_0$ generates, in the sense that each element of the partition $\vee_{n = 0}^\infty T_0^{-n}\sP_0$ is a singleton.
It follows that for every pair of points $x, y \in \Delta_0$ there is a non negative integer $s$, such that $T_0^s(x)$ and $T_0^s(y)$ belong to different elements of the partition $\sP_0$.
We denote by $s_0(x, y)$ the least of such integers~$s$ and call it {\it the separation time} of $x$ and $y$.

We assume furthermore that for each element $\Delta'$ of $\sP_0$ the map $(T_0|_{\Delta'})^{-1}$ is measurable and that the Jacobian $\Jac(T_0)$ of $T_0$ is well defined and positive on a set of full measure of $\Delta_0$.
Moreover, we require that there are constants $C > 0$ and $\beta \in (0, 1)$, such that for almost every $x, y \in \Delta_0$ that belong to the same element of $\sP_0$, we have
$$
\left| \Jac(T_0)(x)/\Jac(T_0)(y) - 1 \right| \le C \beta^{s_0(T_0(x), T_0(y))}.
$$

Let $R$ be a measurable function defined on $\Delta_0$, taking positive integer values and that is constant on each element of $\sP_0$.
Moreover we assume that the greatest common divisor of the values of $R$ is equal to~$1$.
Then put,
\begin{equation*}
  \label{eq:tower}
  \Delta = \{ (z, n) \in \Delta_0 \times \{ 0, 1, \ldots \} \mid n < R(z) \},
\end{equation*}
and endow $\Delta$ with the measure $\measure$, such that for each $n = 0, 1, \ldots$ its restriction to $\Delta_n \= \{ (z, n) \mid z \in \Delta_0, (z, n) \in \Delta \}$ is equal to the pull-back of $\measure_0$ by the map $(z, n) \mapsto z$.
Moreover we define the map $T : \Delta \to \Delta$ by
\begin{eqnarray*}
T(z, n) =
\begin{cases}
(z, n + 1) & \text{ if $n  + 1 < R(z)$} \\
(T_0(z), 0)  & \text{ if $n  + 1 = R(z)$}.
\end{cases}
\end{eqnarray*}

Recall that exponential mixing and Central Limit Theorem were defined in~\S\ref{si:invariant measures}.
\begin{theo }[L.-S. Young~\cite{You}]
With the previous considerations, the following properties hold.
\begin{enumerate}
\item[1.]
If $\int R \, d\measure_0 < + \infty$, then the map $T : \Delta \to \Delta$ admits an invariant probability measure $\rho$ that is that is absolutely continuous with respect to $\measure$.
Moreover, the measure~$\rho$ is ergodic, mixing, and its density with respect to $\measure$ is almost everywhere bounded from below by a positive constant.
\item[2.]
If $\measure_0(\{ z \in \Delta_0 \mid R(z) > m \})$ decreases exponentially fast with~$m$, then the measure $\rho$ is exponentially mixing and the \CLT{} holds for~$\rho$.
\end{enumerate}
\end{theo }

We will also need the following general lemma, whose proof is below.
\begin{lemm}\label{l:induced IFS}
Let $f$ be a rational map and let $(\hV, V)$ be a nice couple for $f$ such that the corresponding canonical induced map $F: D \to V$ satisfies the conclusions of the Key Lemma.
Given $\tc \in \sC$, let $\tF : \tD \to V^{\tc}$ be the first return map of $F$ to $V^{\tc}$ and denote by $J(\tF)$ the maximal invariant set of $\tF$.
Then we have,
$$
\HD((J(f) \cap V^{\tc}) \setminus J(\tF)) < \HD(J(f)),
$$
$\HD(J(f) \setminus \cup_{n \ge 0} f^{-n}(J(\tF))) < \HD(J(f))$ and there is $\talpha \in (0, \HD(J(f)))$ such that,
\begin{equation}\label{e:finite pressure at talpha}
 \sum_{W \text{ c.c. of } \tD} \diam(W)^{\talpha} < + \infty. 
\end{equation}
\end{lemm}
To prove theorems~\ref{t:invariant measures} and~\ref{t:statistical properties}, let~$f$ be a rational map satisfying the \ExpShrink{} condition and let~$\mu$ be the conformal probability measure of exponent $\alpha(f) = \HD(J(f))$ for~$f$, given by Theorem~\ref{t:conformal measures}, so that $\mu$ is supported on the conical Julia set of~$f$ and $\HD(\mu) = \HD(J(f))$.
We will use several times that the measure $\mu$ does not charge sets whose Hausdorff dimension is strictly less than $\HD(J(f))$.
To see this, suppose by contradiction that there is a set $X$ such that $\mu(X) > 0$ and $\HD(X) < \HD(J(f))$.
Then the set $Y = \cup_{n \ge 0} f^n(X)$ has full measure with respect to~$\mu$ (Proposition~\ref{p:conical conformal}) and we have that $\HD(Y) = \HD(X) < \HD(J(f))$.
But this contradicts $\HD(\mu) = \HD(J(f))$.

Let $(\hV, V)$ be a nice couple for~$f$ satisfying the conclusions of Lemma~\ref{l:mixingness}, Theorem~\ref{t:induced map} and of the Key Lemma, as in the proof of Theorem~\ref{t:conformal measures} in~\S\ref{ss:proof of A}.
Let $F : D \to V$ be the canonical induced map associated to $(\hV, V)$.
Moreover, let $\tc \in \sC$ be given by Lemma~\ref{l:mixingness} and let $\tF : \tD \to V^{\tc}$ be the first return map~$F$ to $V^{\tc}$.
We denote by $R$ the return time function of $\tF$ with respect to~$f$, so that $\tF \equiv f^R$ on $\tD$.

Put $\Delta_0 = J(\tF)$, $\measure_0 = \mu|_{J(\tF)}$, $T_0 = \tF|_{J(\tF)}$ and $\sP_0 = \{ W \cap J(\tF) \mid W \text{ c.c. of } \tD \}$.
Notice that the measure $\measure_0$ is non-zero, because $\HD((J(F) \cap V^{\tc}) \setminus J(\tF)) < \HD(J(f))$ (Lemma~\ref{l:induced IFS}) and hence $\mu(J(\tF)) = \mu(J(F) \cap V^c)  > 0$.
We will now verify that all the hypothesis of Young's theorem are verified for this choice of $\Delta_0$, $T_0$, $\sP_0$ and~$R$.

First notice that $\tF$ is an induced map of $f$ in the sense of Appendix~\ref{a:induced maps}.
Hence, the bounded distortion property and the fact that the partition $\sP_0$ is generating for $T_0 = \tF|_{J(\tF)}$, follow from~\S\ref{sa:generalities of induced}.
On the other hand, the final statement of Lemma~\ref{l:mixingness} implies that the greatest common divisor of the values of~$R$ is equal to~$1$.
It remains to prove the ``tail estimate'' in part~$2$ of Young's theorem.
By the bounded distortion property of $\tF$, it follows that there is a constant $C_0 > 1$ such that for every connected component $W$ of $\tD$ we have
$$
C_0^{-1} \diam(W)^{\alpha(f)} \le \mu(W) \le C_0 \diam(W)^{\alpha(f)}.
$$
Note that for each connected component $W$ of $\tD$ and for every $z \in W$ we have $R(z) = m_W$.
Moreover, since for each $z \in D$ we have $|F'(z)| \ge \la^{m(z)}$ it follows that for every connected component $W$ of $\tD$ we have that $\diam(W) \le \la^{-m_W} \diam(V^{c(W)})$.
As $\talpha < \HD(J(f)) = \alpha(f)$, for each positive integer $m$ we have
\begin{equation*}
\begin{split}
\mu(\{z \in \tD \mid R(z) \ge m \})
& =
\sum_{W \text{ c.c. } \tD, m_W \ge m} \mu(W)
\\ & \le
C_0 \sum_{W \text{ c.c. } \tD, m_W \ge m} \diam(W)^{\alpha(f)}
\\ & \le
C_1 \sum_{W \text{ c.c. } \tD, m_W \ge m} (\lambda^{- (\alpha(f) - \talpha)})^{m_W} \diam(W)^{\talpha}
\\ & \le
C_1 (\lambda^{- (\alpha(f) - \talpha)})^m \sum_{W \text{ c.c. } \tD} \diam(W)^{\talpha},
\end{split}
\end{equation*}
where $C_1 \= C_0 \left(\max_{c \in \sCJ} \diam(V^c)\right)^{\alpha(f) - \talpha}$.
By~\eqref{e:finite pressure at talpha} it follows that there is a constant $C_2 > 0$ such that, letting $\theta = \lambda^{- (\alpha(f) - \talpha)} \in (0, 1)$, we have
$$
\mu(\{z \in D \mid R(z) \ge m \})
\le C_2 \theta^m.
$$

Thus we have verified all the hypothesis of Young's theorem.
Let $\rho$ be the invariant measure given by Young's theorem and consider the projection $\pi : \Delta \to \CC$ defined by $\pi(z, n) = f^n(z)$.
We have $f \circ \pi = \pi \circ T$ and therefore the measure $\pi_*\rho$ is invariant by $f$, it is exponentially mixing and the \CLT{} holds for this measure.
It remains to show that this measure is absolutely continuous with respect to $\mu$ and that its density is almost everywhere bounded from below by a positive constant.
First notice that by definition of $\measure$ we have $\pi_* \measure|_{J(\tF) \times \{ 0 \}} = \mu|_{J(\tF)}$ and that for every connected component $W$ of $\tD$ and every $n \in \{ 0, \ldots, m_W - 1 \}$ the measure
$$
\mu_{W, n}
\=
\pi_* \measure|_{(W \cap J(\tF)) \times \{ n \}}
=
(f^n)_* \mu|_{W \cap J(\tF)}
$$
is absolutely continuous with respect to $\mu$, with density $J_{W,n}:=|(f^n)'|^{-\alpha}$ on $f^n(W)$, and~$0$ in the rest of $\CC$.
When we sum these measures over all possible $W$ and~$n$, we obtain a series of measures $\sum_{W,n} \mu_{W,n}$, that converges to $\pi_* \measure$ as linear functionals on continuous functions defined on $J(f)$.
If $u : J(f) \to \R$ is the constant function equal to~$1$, then by the Monotone Convergence Theorem we have,
$$
+ \infty > \pi_*\measure (J(f)) = \sum_{W,n} \int J_{W,n} \,d\mu =
\int \sum_{W,n} J_{W,n}\,d\mu.
$$
Hence the function $\sum_{W,n} J_{W,n}$ is $\mu$-integrable.
Using the Monotone Convergence Theorem again we have that for every continuous function $u : J(f) \to \R$,
\begin{multline*}
\int u d\pi_*\measure
=
\sum_{W,n} \int u d\mu_{W,n}
=
\sum_{W,n} \int uJ_{W,n}\, d\mu
= \\ =
\int \sum_{W,n}  uJ_{W,n}\, d\mu
=
\int u(\sum_{W,n} J_{W,n})\, d\mu.
\end{multline*}
Thus, it follows that $\pi_* \measure$ is absolutely continuous with respect to $\mu$, with density $\sum_{W,n} J_{W,n}$.

As $\rho$ is absolutely continuous with respect to $\measure$, we have that $\pi_*\rho$ is absolutely continuous with respect to $\mu$.
Let $h$ be the density of $\pi_* \rho$ with respect to~$\mu$.
Since $\pi_0 \measure|_{J(\tF) \times \{ 0 \}} = \mu|_{J(\tF)}$, it follows by Young's theorem that there is a constant $c > 0$, such that we have $h > c$ almost everywhere on $J(\tF)$.
As $\mu(V^{\tc} \setminus J(\tF)) = 0$ (because the Hausdorff dimension of this set is strictly less than that of $J(f)$), it follows that we have $h > c$ almost everywhere on $V^{\tc}$.
Let $N$ be a positive integer such that $J(f) \subset f^N(V^{\tc})$.
By the invariance of $\mu$, for almost every $z \in J(f)$ we have
$$
h(z)
=
\sum_{f^N(y) = z} h(y)|(f^N)'(y)|^{-\alpha(f)}
\ge c \left( \sup_{\CC} |(f^N)'| \right)^{-\alpha(f)}.
$$
This finishes the proof of theorems~\ref{t:invariant measures} and~\ref{t:statistical properties}. 
\begin{rema}\label{r:Markov versus IFS}
  When $\sC$ contains exactly one element Lemma~\ref{l:induced IFS} is not necessary, because in this case $\tF = F$.
When $\sC$ contains more than one element we cannot apply Young's theorem directly to $\Delta_0 = J(F)$, $T_0 = F$ and $\sP_0 = \{ W \cap J(F) \mid W \text{ c.c. of } D \}$, because in this case the image by $F$ of an element of $\sP_0$ is of the form $V^c \cap \Delta_0$, for some $c \in \sC$, and it is not equal to~$\Delta_0$.
However Young's theorem extends to this more general setting, as shown in~\cite[Th\'eor\`eme~$2.3.6$ and Remarque $2.3.7$]{G} or \cite[Theorem~$1.1$]{BM}, and we can also apply this more general result directly to~$F$.
\end{rema}
\begin{proof}[Proof of Lemma~\ref{l:induced IFS}.]
  Let $G$ be the restriction of $F$ to $(D \setminus V^{\tc}) \setminus F^{-1}(V^{\tc})$ and notice that $G$ is an induced map of~$f$ in the sense of Appendix~\ref{a:induced maps}.
As~$F$ is topologically mixing (Lemma~\ref{l:mixingness}), the domain of $G$ is strictly smaller than $D \setminus V^{\tc}$ and therefore we have that $\HD(J(G)) < \HD(J(F))$~\cite[Theorem~4.7, p. 134]{MU}.
It is easy to see that the set $(J(F) \cap V^{\tc}) \setminus J(\tF)$ is equal to the preimage of $J(G)$ by $F|_{V^{\tc}}$.
As each inverse branch of $F$ is Lipschitz, it follows that
$$
\HD((J(F) \cap V^{\tc}) \setminus J(\tF)) \le \HD(J(G)) < \HD(J(F)) = \HD(J(f))
$$
and that $\HD((J(f) \cap V^{\tc}) < \HD(J(f))$.
As $\HD(J(f) \cap K(V)) < \HD(J(f))$ (Lemma~\ref{l:hd}), it also follows that $\HD(J(f) \setminus \cup_{n \ge 0} f^{-n}(J(\tF))) < \HD(J(f))$.

As by hypothesis $F$ satisfies the conclusions of the Key Lemma, it follows that the pressure function of~$F$ is finite at $\alpha$ and that it vanishes at $\HD(J(F)) = \HD(J(f))$.
Thus, the pressure function of $G$ is finite at $\alpha$, and from the inequality $\HD(J(G)) < \HD(J(F))$ it follows that the pressure function of $G$ is negative at $\HD(J(F))$.
Thus there is $\talpha \in (\alpha, \HD(J(f)))$ where the pressure function of~$G$ is negative.

Let $D'$ be the subset of $D$ of all those points~$z$ for which there is a positive integer~$m$ such that $F^m$ is defined at~$z$ and such that $F^m(z) \in V^{\tc}$.
For each $z \in D'$ we denote by $m'(z)$ the least value of such~$m$, so that $\tF \equiv F^{m'}$ on~$\tD$.
Note that $D' \cap V^{\tc} = \tD$.
The function $m'$ is constant on each connected component of $D'$.
For a connected component $W$ of $D'$ we denote by $m_W'$ the common value of $m'$ on~$W$.
If $W$ is a connected component of $D'$ such that $m_W' = 1$, then $F(W) = V^{\tc}$.
In the case $m_W' > 1$, the set $F(W)$ is a connected component of $D'$ and $m_{F(W)}' = m_W' - 1$.
For each $c \in \sC \setminus \{ \tc \}$ fix a point $z_c \in V^c$.
Thus, if we denote by $C_0 > 0$ the distortion constant of~$F$, then we have
\begin{multline*}
  \sum_{W \text{ c.c. of } D' \setminus V^{\tc}} \diam(W)^{\talpha}
  \le C_0^{\talpha} \left(\sum_{W \text{ c.c. of } D' \setminus V^{\tc}, m_W' = 1} \diam(W)^{\talpha} \right)
  \cdot \\ \cdot
\left( \sum_{m \ge 1} \sum_{c \in \cC \setminus \{ \tc \}} \sum_{y \in (G)^{-m}(z_c)} |(F^m)'(y)|^{-\talpha} \right)
\end{multline*}
Note that the first factor is finite by part~$2$ of the Key Lemma and that the second factor is finite because the pressure function of $G$ is negative at~$\talpha$.
Now, if $W'$ is a connected component of $D$ that is contained in $V^{\tc}$ and whose image by $F$ is not equal to $V^{\tc}$, then we have $F(W' \cap \tD) = D' \cap V^{c(W')}$ and the sum,
$$
\sum_{W \text{ c.c. of } W' \cap \tD} \diam(W)^{\talpha}
$$
is less than a distortion constant times,
$$
\diam(W')^{\talpha} \cdot \left( \sum_{W \text{ c.c. of } D' \cap V^{c(W')}} \diam(W)^{\talpha} \right).
$$
Then the estimate~\eqref{e:finite pressure at talpha} follows from part~$2$ of the Key Lemma.
\end{proof}
\subsection{Proof of Theorem~\ref{t:reverse implication}.}\label{ss:proof of D}
Let $f$ be a rational map having an exponentially mixing invariant measure $\nu$, that is absolutely continuous with respect to a conformal invariant measure $\mu$ of~$f$.
Moreover we assume that there is a constant $c > 0$ such that the density of~$\nu$ with respect to~$\mu$ is almost everywhere larger than~$c$.

We will show that there is a constant $\lambda > 1$, such that for every positive integer~$n$ and every repelling periodic point $p$ of~$f$ of period~$n$ we have $|(f^n)'(p)| \ge \lambda^n$.
By~\cite{PRS} this last property is equivalent to the \TCEC, so this proves the theorem.

As $\nu$ is exponentially mixing, there are constants $C > 0$ and $\rho \in (0, 1)$ such that for every continuous function $\varphi : J(f) \to \R$, every Lipschitz function $\psi : J(f) \to \R$ and every positive integer~$n$ we have~\ref{e:correlation bound}.
Let $n$ be a positive integer and let $p$ be a repelling periodic point of~$f$ of period~$n$.
Let $\phi$ be a local inverse of $f^n$ at~$p$ that fixes~$p$.
Let $r > 0$ be sufficiently small such that $\phi$ is defined on the ball $B(p, r)$ and such that $\phi(B(p, r)) \subset B(p, r)$.

By the Fatou-Sullivan classification of connected components of the Fatou set~\cite{Bea,CG,Mi}, is easy to see that $\nu$, and hence $\mu$, is supported on the Julia set of $f$.
Therefore the topological support of $\mu$, and hence of $\nu$, is equal to the Julia set of $f$.
It follows that there is $\e > 0$ such that $\nu(B(p, r) \setminus B(p, \e r)) > 0$.
Thus there is a bounded and measurable function $\varphi$ that is constant equal to~$1$ on $B(p, \e r)$, that is constant equal to~$0$ outside of $B(p, r)$ and such that $\int \varphi d\nu = 0$.

By Koebe Distortion Theorem there is a constant $\eta > 1$ such that for every positive integer $k$ there is $r' > 0$ such that,
$$
\eta^{-1} r |(f^{kn})'(p)|^{-1} < r' < \eta r |(f^{kn})'(p)|^{-1}),
$$
and
$$
B(p, r') \subset \phi^k(B(p, \e r)) \subset B(p, \eta r').
$$
By the conformality of $\mu$ it also follows that, if we denote by $\alpha$ the exponent of~$\mu$, we can take $\eta > 1$ large enough so that $\mu(B(p, \tfrac{1}{2} r')) \ge \eta^{-1} |(f^{kn})'(p)|^{-\alpha}$.
For such $k$ and~$r'$, let $\psi : J(f) \to \R$ be a Lipschitz function that is equal to~$0$ outside $B(p, r')$, that is positive on $B(p, r')$, and that is equal to~$1$ on $B(p, \tfrac{1}{2}r')$.
Moreover, we can take such a $\varphi$ in such way that for some constant $C' > 0$ independent of $k$ we have $\| \psi \|_{\Lip} \le C' |(f^{kn})'(p)|$.
It follows that
\begin{multline*}
\left| \int (\varphi \circ f^n) \cdot \psi d\nu - \int \varphi d\nu \int \psi d\nu \right|
=
 \left| \int (\varphi \circ f^n) \cdot \psi d\nu \right|
\ge \\ \ge
\nu(B(p, \tfrac{1}{2} r'))
\ge
c \mu(B(p, \tfrac{1}{2}r'))
\ge
c \eta^{-1} |(f^{kn})'(p)|^{-\alpha}.
\end{multline*}
Then, the inequality~\eqref{e:correlation bound} implies that,
$$
|(f^{kn})'(p)|^{-\alpha}
\le 
c^{-1} \eta C \cdot \left( \sup_{J(f)} |\varphi| \right) \cdot \| \psi \|_{\Lip} \cdot \rho^n
\le
C'' |(f^{kn})'(p)| \rho^{kn},
$$
where $C'' = c^{-1} \eta C \cdot \left( \sup_{J(f)} |\varphi| \right) C'$.
Thus we have,
$$
|(f^n)'(p)|^k
=
|(f^{kn})'(p)|
\ge
C'' \left( \rho^{-\tfrac{1}{1 + \alpha}} \right)^{kn}.
$$
As this holds for every positive integer $k$, it follows that $|(f^n)'(p)| \ge \left( \rho^{-\tfrac{1}{1 + \alpha}} \right)^n$.
This shows the desired assertion with $\lambda \= \rho^{-\tfrac{1}{1 + \alpha}}$.
\subsection{Characterizations of the invariant measure.}\label{ss:characterizations}
The purpose of this subsection is to prove the following proposition.
\begin{prop}\label{p:HD of invariant}
  For a rational map~$f$ satisfying the \TCEC, the measure given by Theorem~\ref{t:invariant measures} is the unique invariant measure supported on $J(f)$ whose Hausdorff dimension is equal to $\HD(J(f))$.
\end{prop}
Before the proof of this proposition we consider the following corollary.
For an invariant measure~$\nu$ we will denote by $h_\nu$ its metric entropy and by $\chi_\nu \= \int \ln |f'| d\nu$ its Lyapunov exponent.
An $f$-invariant probability measure is called an \textit{equilibrium state with potential $- \HD(J(f)) \ln |f'|$ of~$f$}, if it is supported on $J(f)$ and if it maximizes
\begin{equation}\label{e:pressure}
h_\nu + \int - \HD(J(f)) \ln |f'| d \nu
= h_\nu - \HD(J(f)) \chi_\nu,
\end{equation}
among all $f$-invariant probability measures $\nu$ that are supported on~$J(f)$.
\begin{coro}\label{c:invariant as equilibrium}
  For a rational map~$f$ satisfying the \TCEC, the measure given by Theorem~\ref{t:invariant measures} is the unique equilibrium state with potential $-\HD(J(f)) \ln |f'|$ of~$f$. 
\end{coro}
\begin{proof}  
When $f$ satisfies the \TCEC, the supremum of~\eqref{e:pressure} is equal to~$0$~\cite{Prconical,PRSpressure}. 
Thus, in this case an invariant probability measure $\nu$ is an equilibrium state with potential $- \HD(J(f)) \ln |f'|$ of~$f$, if and only if $\nu$ is supported on~$J(f)$ and if $h_\nu = \HD(J(f)) \cdot \chi_\nu$.

Let $\nu$ be an $f$-invariant measure supported on~$J(f)$, and consider its ergodic decomposition $\nu = \int \nu_\xi d \tnu(\xi)$.
As~$f$ satisfies the \TCEC, for each $\xi$ the Lyapunov exponent of~$\nu_\xi$ is positive~\cite{PRS} and therefore $h_{\nu_\xi}/\chi_{\nu_\xi} = \HD(\nu_\xi)$~\cite{Mane,PUbook}. 
Thus we have
$$
h_\nu = \int \HD(\nu_\xi) \chi_{\nu_\xi} d \tnu(\xi),
$$
and we conclude that an invariant measure is an equilibrium state of~$f$ with potential $- \HD(J(f)) \ln |f'|$, if and only if the Hausdorff dimension of almost every ergodic component is equal to $\HD(J(f))$.

As the conformal probability measure of minimal exponent of~$f$ has Hausdorff dimension equal to~$\HD(J(f))$ (Theorem~\ref{t:conformal measures}), it follows that the Hausdorff dimension of the measure given by Theorem~\ref{t:invariant measures} is equal to $\HD(J(f))$.
As this measure is ergodic, it follows that it is an equilibrium state with potential $- \HD(J(f)) \ln |f'|$ of~$f$.
The uniqueness follows from Proposition~\ref{p:HD of invariant}.
\end{proof}
\begin{proof}[Proof of Proposition~\ref{p:HD of invariant}.]
We keep the considerations and notation of the proof of theorems~\ref{t:invariant measures} and~\ref{t:statistical properties}.
Recall that the measure given by Theorem~\ref{t:invariant measures} was obtained as the projection by~$\pi$ of a $T$-invariant measure supported on $\Delta$.
The restriction of this measure to $J(\tF) \times \{ 0 \} \sim J(\tF)$ is a $\tF$-invariant measure whose Hausdorff dimension is equal to $\HD(J(f))$.
By Theorem~\ref{t:conformal for induced} in Appendix~\ref{a:induced maps}, there is a unique $\tF$-invariant probability measure whose Hausdorff dimension is equal to $\HD(J(\tF)) = \HD(J(f))$.
Thus, we just need to show that if $\nu$ is an ergodic $f$-invariant measure supported on the Julia set of $J(f)$ whose Hausdorff dimension is equal to $\HD(J(f))$, then $\nu$ is obtained from an $\tF$-invariant measure whose Hausdorff dimension equal to $\HD(J(f))$, in the same way as it was described above.

If $\nu$ is an ergodic $f$ invariant measure whose Hausdorff dimension is equal to $\HD(J(f))$, then Lemma~\ref{l:induced IFS} implies that $\nu$ is supported on $\cup_{n \ge 0} f^{-n}(J(\tF))$.
Thus we have $\nu(J(\tF)) > 0$ and by~\cite{Zw} there is an $\tF$-invariant measure $\rho_0$, which is absolutely continuous with respect to~$\nu$ and such that~$\nu$ is obtained by first extending $\rho_0$ to a $T$-invariant measure on $\Delta$, and then by projecting it by $\pi$.
As $\rho_0$ is absolutely continuous with respect to $\nu$ it follows that the Hausdorff dimension of~$\rho_0$ is equal to $\HD(J(f)) = \HD(J(\tF))$.
This completes the proof.
\end{proof}
\appendix
\section{Induced maps.}\label{a:induced maps}
In this appendix we study a class of induced maps of a given rational map.
We show in particular that these maps fall into the category of maps studied in~\cite{MUbook}, and gather in Theorem~\ref{t:conformal for induced} several results of this book.

Throughout all this section we fix a rational map~$f$.
\subsection{Definition and general properties of induced maps.}\label{sa:generalities of induced}
Recall that for a nice set $V = \cup_{c \in \sC} V^c$ for~$f$ and for each pull-back~$W$ of~$V$, we denote by $c(W) \in \sC$ the critical point and by $m_W \ge 0$ the integer such that $f^{m_W}(W) = V^{c(W)}$, see~\S\ref{ss:nice sets}.
Moreover, if $(\hV, V)$ is a nice couple for~$f$, then for each pull-back $W$ of~$V$ we denote by~$\hW$ the unique pull-back of~$\hV$ that contains $W$ and such that $m_{\hW} = m_W$.
\begin{defi}\label{d:induced map}
Let $f$ be a rational map and let $(\hV, V)$ be a nice couple for~$f$.
We will say that a map $F : D \to V$, with $D \subset V$, is \textsf{\textit{induced by}} $f$, if the following properties hold: Each connected component $W$ of $D$ is a pull-back of $V$ and we have $F|_W = f^{m_W}|_W$ and $f^{m_W}$ is univalent on~$\hW$.
\end{defi}
It is straightforward to check that every induced map satisfies the following properties.

\begin{enumerate}
\item[\axiomM{1}]{\it Markov property}.
For every connected component $W$ of $D$, the map $F|_W$ is a biholomorphism between $W$ and $V^{c(W)}$.
\item[\axiomM{2}]{\it Univalent extension}.
For each $c \in \sCJ$ and for every connected component $W$ of $D$ satisfying $c(W) = c$, the inverse of $F|_W$ extends to a biholomorphism between $\hV^c$ and $\hW$.
\item[\axiomM{3}]{\it Strong Separation}.
The connected components of~$D$ have pairwise disjoint closures.
\end{enumerate}

Denote by $\fD$ the collection of connected components of $D$ and for each $c \in \sCJ$ denote by $\fD^c$ the collection of all elements of $\fD$ contained in $V^c$, so that $\fD = \sqcup_{c \in \sCJ} \fD^c$.
A word on the alphabet $\fD$ will be called {\it admissible} if for every pair of consecutive letters $W, W' \in \fD$ we have $W \in \fD^{c(W')}$.
For a given integer $n \ge 1$ we denote by~$E^n$ the collection of all admissible words of length~$n$ and we set $E^* = \sqcup_{n \ge 1} E^n$.
Moreover we denote by $E^\infty$ the collection of all infinite admissible words of the form $W_1W_2 \ldots$~.
Given an integer $n \ge 1$ and an infinite word $\underline{W} = W_1 W_2 \ldots \in E^\infty$ or a finite word $\underline{W} = W_1 \ldots W_m \in E^*$ of length $m \ge n$, put $\underline{W}|_n = W_1 \ldots W_n$.

Given $W \in \fD$, denote by $\phi_W$ the holomorphic extension to $\hV^{c(W)}$, of the inverse of $F|_W$, that is given by property \axiomM{2}.
For a finite word $\underline{W} = W_1 \ldots W_n \in E^*$ put $c(\underline{W}) \= c(W_n)$.
Note that the composition
$$
\phi_{\underline{W}} \= \phi_{W_1} \circ \ldots \circ \phi_{W_n}
$$
is well defined and univalent on $\hV^{c(\uW)}$ and takes images into $V$.
Moreover, put
$$
D_{\uW} \= \phi_{\uW} (V^c), \
\hD_{\uW} \= \phi_{\uW} (\hV^c)
\ \mbox{ and } \
A_{\uW} \= \hD_{\uW} \setminus D_{\uW}.
$$
Note that $A_{\uW}$ is an annulus of the same modulus as the annulus $\hV^{c(\uW)} \setminus \ov{V^{c(\uW)}}$ and that, when $\uW \= W \in \fD$, we have $D_{W} = W$ and $\hD_{W} = \hW$.

\subsubsection*{Bounded distortion.}
For each $\uW \in E^*$ the map $\phi_{\uW}$ is defined and univalent on $\hV^{c(\uW)}$.
It follows by Koebe Distortion Theorem that the distortion of $\phi_{\underline{W}}$ on $V^{c(\uW)}$ is bounded independently of $\underline{W}$.

Given $n \ge 0$ note that the map $F^n$ is well defined on $\sqcup_{\uW \in E^n} D_{\uW}$.
Moreover, for each $\uW \in E^n$ the restriction of $F$ to $D_{\uW}$ is equal the inverse of $\phi_{\underline{W}}$ on this set.
It follows that the distortion of $F^n$ on each $D_{\uW}$ is bounded independently of $n$ and of $\uW$.
\subsubsection*{Expansion.}
For each $c \in \sCJ$ endow $\hV^c$ with the corresponding hyperbolic metric.
The restriction of this metric to $V^c$ is comparable to the spherical metric on $\CC$.
By Schwarz-Pick Lemma it follows that for each $W \in \fD$ the holomorphic map $\phi_W : \hV^{c(W)} \to V$ decreases the hyperbolic metric by a factor $s \in (0, 1)$, independent of $c$ and $W$.
So there is a constant $C_D > 0$ such that for every finite word $\underline{W} \in E^n$ the spherical diameter of $D_{\underline{W}}$ is at most $C_D \cdot s^n$. 
\subsubsection*{Maximal invariant set.}
For every infinite word $\underline{W} \in E^\infty$ and every $n \ge 1$, we have
$$
\ov{D}_{\uW|_{n + 1}} \subset D_{\uW|_n}
\ \mbox{ and } \
\diam(D_{\uW|_n}) \le C_D \cdot s^n.
$$
It follows that the intersection \ $\cap_{n \ge 1} D_{\underline{W}|_n}$ \ is a singleton.
We denote the corresponding point by $\pi(\underline{W})$.
So $\pi$ defines a map from $E^\infty$ to $\ov{\C}$; we denote by $J(F) \= \pi(E^\infty)$ the image of $\pi$.
Note that the set $J(F)$ is equal to the maximal invariant set of $F$.

It follows from condition~\axiomM{3} that for every integer $n \ge 1$ and for distinct $W, W' \in E^n$, the closure of the sets $D_{\underline{W}}$ and $D_{\underline{W}'}$ are disjoint.
Therefore~$\pi$ induces a bijection between $E^\infty$ and $J(F)$.

The following technical property is important to use the results of~\cite{MUbook}.
\begin{prop}\label{p:induced maps}
Every induced map satisfies the following property
\begin{enumerate}
\item[\axiomM{4}]
There is a constant $C_M > 0$ such that for every $\kappa \in (0, 1)$ and every ball $B$ of $\CC$, the following property holds.
Every collection of pairwise disjoint sets of the form $D_{\underline{W}}$, with $\underline{W} \in E^*$, intersecting $B$ and with diameter at least $\kappa \cdot \diam(B)$, has cardinality at most $C_M \kappa^{-2}$.   
\end{enumerate}
\end{prop}
The proof of this proposition is based on the following lemma.
\begin{lemm}\label{l:distribution}
Let $f$ by a rational map and let $F : D \to V$ be a map induced by $f$.
Then the following properties hold.
\begin{enumerate}
\item[1.]
There exists a constant $C_0 > 0$ such that if $\uW$ and $\uW' \in E^*$ are such that $D_{\uW}$ and $D_{\uW'}$ are disjoint, then
$$
\dist(D_{\uW'}, D_{\uW}) \le C_0 \cdot \diam(D_{\uW}),
\ \mbox{ implies } \
\hD_{\uW'} \subset A_{\uW}.
$$
\item[2.]
Let $B$ be a ball.
Then for every pair $\uW, \uW' \in E^*$ such that the sets $D_{\uW}$ and $D_{\uW'}$ are disjoint and intersect $B$, we have either
$$
\diam(D_{\uW}) < C_0^{-1} \diam(B)
\ \mbox{ or } \
\diam(D_{\uW'}) < C_0^{-1} \diam(B).
$$
\end{enumerate}
\end{lemm}
\begin{proof}
\

\partn{1}
Observe that for every $\uW \in E^*$ the set $A_{\uW}$ is an annulus whose modulus is at least
$$
\min_{c \in \sCJ} \mod (\hV^c \setminus \ov{V^c}) > 0.
$$
It follows that there is a constant $C_0 > 0$ such that every point whose distance to $D_{\uW}$ is at most $C_0 \cdot \diam(D_{\uW})$ is contained in $\widehat{D}_{\uW}$.
So, if $\uW'$ is an element of $E^*$ such that $D_{\uW'}$ is disjoint from $D_{\uW}$ and at distance at most $C_0 \cdot \diam(D_{\uW})$ from $D_{\uW}$, then $D_{\uW'}$ intersects $\widehat{D}_{\uW}$.
As the sets $\widehat{D}_{\uW'}$ and $\widehat{D}_{\uW}$ are both pull-backs of~$\hV$, it follows that in this case we have $\widehat{D}_{\uW'} \subset \widehat{D}_{\uW} \setminus \ov{D}_{\uW} = A_{\uW}$.

\partn{2}
Assume by contradiction that there are such $\uW$ and $\uW'$ in $E^*$ for which $\diam(D_{\uW})$ and $\diam(D_{\uW'})$ are both bigger than or equal to $C_0^{-1} \diam(B)$. 
Then
$$
\dist(D_{\uW}, D_{\uW'})
\le
\diam(B)
\le
C_0 \cdot \min \left\{ \diam(D_{\uW}), \diam(D_{\uW'}) \right\},
$$
and part~$1$ implies that $\widehat{D}_{\uW} \subset A_{\uW'}$ and that $\widehat{D}_{\uW'} \subset A_{\uW}$.
So we get a contradiction that proves the assertion.
\end{proof}

\begin{proof}[Proof of Proposition~\ref{p:induced maps}.]
By the bounded distortion property there is a constant $C > 0$ such that for every $\uW \in E^*$ we have
$$
\Area(D_{\uW}) \ge C^{-1} \diam(D_{\uW})^2.
$$
Taking $C$ larger if necessary we assume that for every ball $B$ in $\ov{\C}$ we have $\Area(B) \le C \diam(B)^2$.

To prove property \axiomM{4} let $\kappa \in (0, 1)$ and let $B$ be a ball.
Moreover let $\fF$ be a finite collection of elements of $E^*$, such that the sets $D_{\uW}$, for $\uW \in \fF$, are pairwise disjoint and such that for every $\uW \in \fF$ we have $\diam(D_{\uW}) \ge \kappa \cdot \diam(B)$.
It follows that for every $\uW \in \fF$ we have $\Area(D_{\uW}) \ge C^{-1} \kappa^2 \cdot \diam(D_{\uW})$.
By part~$2$ of Lemma~\ref{l:distribution} it follows that there is at most one element $\uW$ of $\fF$ such that $\diam(D_{\uW}) \ge C_0^{-1} \diam(B)$.
Let $\tB$ be the ball with the same center as $B$ and with radius equal to $C_0^{-1} \diam(B)$.
So for every $\uW \in \fF$ we have $D_{\uW} \subset \tB$, with at most one exception.
Therefore
\begin{multline*}
\# \fF
\le
1 + \Area(\tB) / (C^{-1} \kappa^2 \cdot \diam(B)^2) \\
 \le
1 + 4C_0^{-2}C^2 \kappa^{-2}
\le 
(1 + 4C_0^{-2} C^2) \kappa^{-2}.
\end{multline*}
\end{proof}
\subsection{Pressure function.}
Fix an induced map $F : D \to V$ of~$f$.
Then for $t \ge 0$ we define
$$
Z_n(\fD, t) \= \sum_{\underline{W} \in E^n} \left( \sup \left\{ |\phi_{\underline{W}}'(z)| \mid z \in V^{c(\uW)} \right\} \right)^t.
$$
The proof of the following lemma is standard, see for example lemmas~$2.1.1$ and~$2.1.2$ of~\cite{MU}.
\begin{lemm}
For every $t \ge 0$ the sequence $( \ln Z_n(\fD, t) )_{n \ge 1}$ is sub-additive and
$$
P(t)
:=
\lim_{n \to + \infty} \tfrac{1}{n} \ln Z_n(\fD, t)
=
\inf \left\{ \tfrac{1}{n} \ln Z_n(\fD, t) \mid  n \ge 1 \right\}.
$$
\end{lemm}
The function $P : [0, + \infty) \to [-\infty, +\infty]$ so defined is called the {\it pressure function}.
It is easy to see that for every $t \ge 0$ the sequence $( \tfrac{1}{n} \ln Z_n(\fD, t))_{n \ge 1}$ is uniformly bounded from below.
In particular the function $P(\fD, \cdot)$ does not take the value $-\infty$.
Note however that if $\fD$ is infinite, then $P(0) = + \infty$.

The proof of the following lemma is straightforward.
\begin{lemm}\label{l:regularity of the pressure function}
Put $\theta(F) \= \inf \{ t \ge 0 \mid P(t) < + \infty \}$.
Then the pressure function $P$ is finite, continuous and strictly decreasing to $-\infty$ on $(\theta(F), + \infty)$.
\end{lemm}
It follows from this lemma that there is at most one value of $t \ge 0$ for which the pressure function $P$ vanishes.
Following the terminology of~\cite{MU} we say that $F$ is {\it strongly regular} if there is $t > \theta(F)$ at which $P$ vanishes.
\subsection{Conformal measures.}
Fix an induced map $F : D \to V$ of~$f$.
For a given $t > 0$, a finite Borel measure $\mu$ will be called {\it conformal with exponent $t$} for $F$, if $\mu$ is supported on $J(F)$ and if for every $W \in \fD$ and every Borel set $U$ contained in $W$, we have
$$
\mu(F(U)) = \int_U |F'|^t \ d\mu.
$$
In the following theorem we gather several results of~\cite{MU}, applied to our particular setting.
\begin{theo}[\cite{MU}]\label{t:conformal for induced}
Let $F$ be a topologically mixing and strongly regular induced map and let $h \ge 0$ be the unique zero of the pressure function of~$F$.
Then $h = \HD(J(F))$, there is a unique conformal probability measure $\mu$ of exponent~$h$ for $F$ and this measure satisfies $\HD(\mu) = \HD(J(f))$.
Furthermore, there is a unique invariant probability measure $\rho$ of $F$ that is absolutely continuous with respect to $\mu$ and this measure is also characterized as the unique $F$-invariant probability measure satisfying
$
\HD(\rho) = h.
$
\end{theo}
\begin{proof}
The collection of maps $\Phi = \{ \phi_W \mid W \in \fD \}$ is a {\it Conformal Graph Directed Markov System} ({\it CGDMS} for short), as defined in~$\S 4.2$ of~\cite{MU}, except for the fact that the {\it Cone Condition~(4d)} of~\cite{MU} is replaced here by the weaker condition~\axiomM{4}.
All the results we use from~\cite{MU} only require this weaker condition, as it is explained below.
Moreover note that the hypothesis that $F$ is topologically mixing easily implies that $\Phi$ is finitely primitive in the sense of~\cite{MU}.

The equality $h = \HD(J(F))$ is given by Theorem~$4.2.13$ of~\cite{MU}.
In the proof of this result the Cone Condition is only used to guaranty that the conclusion of Lemma~$4.2.6$ holds.
But the conclusion of this last lemma is our condition~\axiomM{4}.

Theorem~$4.2.9$ of~\cite{MU} implies that there is a conformal measure $\mu$ of exponent~$h$ for~$F$.
Observe that in the proof of Theorem~$4.2.9$ of~\cite{MU} the Cone Condition is only used to prove that the CGDMS in question is conformal-like (Proposition~$4.2.7$ of~\cite{MU}).
But this is an immediate consequence of the Strong Separation Condition~\axiomM{3} (notice that our condition \axiomM{3} is stronger than the {\it Open Set Condition~(4b)} of~\cite{MU}.)

The existence and uniqueness of the absolutely continuous invariant measure is given by theorems~$6.1.2$ and~$6.1.3$ of~\cite{MU}.
Although these results are stated for Conformal Iterated Functions Systems, they apply equally well to CGDMS and to Markov maps.
Finally the equalities $\HD(\mu) = \HD(\rho) = h$ are given by Corollary~$4.4.6$ of~\cite{MU}.
As before, in the proof of these results the Cone Condition is only used to guarantee that the conclusion of Lemma~$4.2.6$ holds.
\end{proof}

\section{Conformal measures via inducing.}\label{a:conformal via inducing}
Fix throughout all this section a rational map $f$ and a nice couple $(\hV, V)$ for~$f$.
In this appendix we study the conformal measures of~$f$, through the canonical induced map $F$ associated to $(\hV, V)$.
In particular we give a sufficient condition on~$F$, for~$f$ to have a conformal measure supported on the conical Julia set.
See~\S\ref{ss:rational conformal measure} and~\S\ref{ss:conical} for the definition of conformal measure and of conical Julia set, respectively.

This appendix is dedicated to the proof of the following result.
See~\S\ref{ss:rational conformal measure} for the definition of~$\alpha(f)$.
\begin{theo}\label{t:conformal via inducing}
Let $f$ be a rational map of degree at least~$2$, let $(\hV, V)$ a nice couple for $f$ and let $F : D \to V$ the canonical induced map associated to $(\hV, V)$.
Assume that~$F$ is topologically mixing and that the following properties hold.
\begin{enumerate}
\item[1.]
For every $c \in \sCJ$ we have $\HD(J(F) \cap V^c) = \alpha(f)$.
\item[2.]
There is $\alpha \in (0, \alpha(f))$ such that
\begin{equation}\label{e:regularity hypothesis}
\sum_{W \text{ c.c. of } D} \diam(W)^{\alpha} < + \infty.
\end{equation}
\end{enumerate}
Then the canonical induced map $F$ is strongly regular in the sense of~\cite{MUbook} and there is a unique conformal probability measure of exponent $\alpha(f)$ for $f$.
Moreover this measure is non-atomic, ergodic, its Hausdorff dimension is equal to $\alpha(f)$ and it is supported on the conical Julia set of~$f$.
\end{theo}
Observe that $J(F)$ is clearly contained in the conical Julia set $\Jcon(f)$ of~$f$.
In~\cite{DU, PrLyapunov} it is shown that $\alpha(f) = \HD(\Jcon(f))$ (see also~\cite{McM}), so by the inclusion $J(F) \subset \Jcon(f)$, we have
$$
\HD(J(F)) \le \HD(\Jcon(f)) = \alpha(f).
$$
So the first hypothesis of the theorem requires, in fact, that for each $c \in \sCJ$ the Hausdorff dimension of $J(F) \cap V^c$ is as large as possible
\subsection{Conformal measures.}\label{ss:rational conformal measure}
For a given $t \ge 0$, we say that a non zero Borel measure $\mu$ is {\it conformal of exponent} $t$ for $f$, if for every Borel subset $U$ of $\CC$ where $f$ is injective, we have
$$
\mu(f(U)) = \int_U |f'|^t d \mu.
$$
By the locally eventually onto property it is easy to see that if the topological support of a conformal measure is contained in the Julia set $J(f)$ of $f$, then the support is in fact equal to $J(f)$.
Note that a conformal measure of exponent $t = 0$ must be supported on the exceptional set of $f$.
So the exponent of a conformal measure supported on $J(f)$ is positive.

It was shown by Sullivan that every rational map admits a conformal measure supported on the Julia set~\cite{Su}.
So the infimum,
\begin{multline*}
\alpha(f)
\=
\inf \{ t > 0 \mid \mbox{ there exists a}
\\
\mbox{conformal measure of exponent $t$ supported on } J(f) \}
\end{multline*}
is well defined and is easy to see that it is realized.
It follows that $\alpha(t)$ is positive, as there is no conformal measure of exponent $t = 0$ supported on the Julia set.

\subsection{The conical Julia set and sub-conformal measures.}\label{ss:conical}
The {\it conical Julia set} of $f$, denoted by $\Jcon(f)$, is by definition the set of all those points $x$ in $J(f)$ for which there exists $\rho(x) > 0$ and arbitrarily large positive integers $n$, such that the pull-back of the ball $B(f^n(x), \rho(x))$ to $x$ by $f^n$ is univalent.
This set is also called {\it radial Julia set}.

We will need the following general result, which is a strengthened version of ~\cite[Theorem~$5.1$]{McM}, \cite[Theorem~$1.2$]{DMNU}, with the same proof.
Given $t \ge 0$ we will say that a Borel measure $\mu$ is {\it sub-conformal of exponent $t$ for} $f$, if for every subset $U$ of $\CC$ on which $f$ is injective we have
\begin{equation}\label{e:sub-conformal}
\int_U |f'|^t d\mu \le \mu(f(U)).
\end{equation}
\begin{prop}\label{p:conical conformal}
If $\mu$ is a sub-conformal measure for $f$ supported on $\Jcon(f)$, whose exponent is at least $\alpha(f)$, then $\mu$ is conformal of exponent $\alpha(f)$ and every other conformal measure of exponent $\alpha(f)$ is proportional to~$\mu$.
Moreover $\mu$ is non-atomic and every subset $X$ of $\CC$ such that $f(X) \subset X$ and $\mu(X) > 0$, has full measure with respect to~$\mu$.
In particular $\mu$ admits at most one absolutely continuous invariant probability measure.
\end{prop}

\subsection{Induced maps and conformal measures.}\label{s:conformal measure for rational map}
For a nice set $V$ for $f$ denote by $R_V$ the first return map to $V$.
\begin{prop}\label{p:spreading}
Let $F$ be the canonical induced map associated to a nice couple $(\hV, V)$ for $f$.
Then every conformal measure of $F$ of exponent greater or equal than $\alpha(f)$ is the restriction to $V$ of a conformal measure of $f$ supported on $\Jcon(f)$.
Moreover these measures have the same Hausdorff dimension.
\end{prop}
The proof of this proposition is provided below, it depends on some lemmas.
Fix a nice couple $(\hV, V)$ for~$f$ and denote by $F : D \to V$ the canonical induced map associated to it.
\begin{lemm}\label{l:spreading}
Denote by $\fD_V$ the collection of connected components of $\CC \setminus K(V)$ and for $W \in \fD_V$ denote by $\phi_W : \hV^{c(W)} \to \hW$ the inverse of $f^{m_W}|_{\hW}$.
Given a conformal measure $\mu$ for $F$ of exponent $t$, for each $W \in \fD_V$ let $\mu_W$ be the measure supported on $W$, defined by
$$
\mu_W(X) = \int_{f^{m_W}(X)} |\phi_W'|^t d\mu.
$$
Then the measure $\sum_{W \in \fD_V} \mu_W$ is supported on $\Jcon(f)$.
If moreover $t \ge \alpha(f)$, then this measure is finite.
\end{lemm}
\begin{proof}
Clearly $J(F) \subset \Jcon(f)$, so $\mu$ and each of the measures $\mu_W$, for $W \in \fD_V$ is supported on $\Jcon(f)$.
To prove the second assertion, let $\hmu$ be a conformal measure for $f$ of exponent $\alpha(f)$.
It follows by Koebe distortion property that for every $W \in \fD_V$ we have $\hmu(W) \sim \diam(W)^{\alpha(f)}$ and $\mu_W(\CC) \sim \diam(W)^t$.
When $t \ge \alpha(f)$, it follows that the measure $\sum_{W \in \fD_V} \mu_W$ is finite.
\end{proof}
\begin{lemm}\label{l:canonical induced map}
The canonical induced map $F$ satisfies the following property:
\begin{enumerate}
\item[\axiomC]
For every connected component $W$ of $D$ we have either $R_V(W) = V^{c(W)}$ or $R_V(W) \subset D$.
\end{enumerate}
\end{lemm}
\begin{proof}
Observe that for every $z \in D$ and every $r = 1, \ldots, m(z) - 1$, we have that $m \= m(z) - r$ is a good time for $f^r(z)$.
So, if $W$ is a connected component of $D$ such that $F \not\equiv R_V$ on $W$, then letting $r$ be such that $R_V$ is equal to $f^r$ on $W$, we have that $m \= m_W - r$ is a good time for every element of $R_V(W)$. So $R_V(W) \subset D$ in this case.
\end{proof}
\begin{proof}[Proof of Proposition~\ref{p:spreading}]
By Lemma~\ref{l:spreading} the measure $\hmu \= \sum_{W \in \fD_V} \mu_W$ is finite and supported on $\Jcon(f)$.
So by Proposition~\ref{p:conical conformal} we just need to prove that $\hmu$ is a sub-conformal measure for $f$.
So, let $U$ be a subset of $\CC$ on which $f$ is injective.
We have to prove that the inequality \eqref{e:sub-conformal} holds.
Clearly the general case follows from the following special cases.

\partn{Case 1} {\sf $U$ is contained in an element $W$ of $\fD_V$, distinct from the $V^c$.}
Then $W' \= f(W) \in \fD_V$, $m_{W'} = m_W - 1$, $c(W') = c(W)$ and $\phi_{W'} = f|_{\hW} \circ \phi_W$.
So
\begin{multline*}
\hmu_{W'}(f(U))
=
\int_{f^{m_{W'}}(f(U))} |\phi'_{W'}|^h \ d \mu
\\ =
\int_{f^{m_W}(U)} |f' \circ \phi_W|^h \cdot |\phi_W'|^h \ d\mu
=
\int_{U} |f'|^h \ d\hmu_W.
\end{multline*}

\partn{Case 2} {\sf $U$ is contained in a connected component $W$ of $D$.}
Let $W'$ be the element of $\fD_V$ containing $f(W)$.
By property \axiomC{} we have $R_V(W) \subset D$.
Let $W''$ be the connected component of $D$ that contains $R_V(W)$.
If $R_V(W) = W''$ then $f$ is equal to $\phi_{W'} \circ F$ on $W$, and then~\eqref{e:sub-conformal} holds with equality.
If $R_V(W)$ is strictly contained in $W''$, then an induction argument using property \axiomC{} shows that there is an integer $n$ such that $F^n$ is well defined on $W''$ and that it coincides with $f^{m_{W''}}$ on this set.
In this case $f$ is equal to $\phi_{W'} \circ (F^n|_{W''})^{-1} \circ F$ on~$W$ and therefore~\eqref{e:sub-conformal} holds with equality.

\partn{Case 3} $U \subset K(V) \cup (V \setminus D)$.
As by definition $\hmu(K(V) \cup (V \setminus D)) = 0$, there is nothing to prove in this case.
\end{proof}
\begin{rema}
The proof of Proposition~\ref{p:spreading} does not give directly that the measure~$\hmu$ is conformal.
In fact, apriori the set $f(V \setminus D)$ might have positive measure with respect to $\hmu$, so in case~$3$ of the proof is not immediate that we have~\eqref{e:sub-conformal} with equality.
\end{rema}
\begin{proof}[Proof of Theorem~\ref{t:conformal via inducing}]
The second assumption implies that the pressure function of~$F$ is finite at~$\alpha$ and, together with the first assumption, this implies that the pressure function of~$F$ is positive at~$\alpha$.
So the induced map $F$ is strongly regular and Theorem~\ref{t:conformal for induced} implies that $F$ admits a conformal measure $\mu$ of exponent $\alpha(f)$ whose Hausdorff dimension is equal to~$\alpha(f)$.
It follows that~$\mu$ is the restriction to $V$ of a conformal measure for~$f$ that is supported on $\Jcon(f)$ and whose Hausdorff dimension is equal to $\alpha(f)$ (Proposition~\ref{p:spreading}).
By Proposition~\ref{p:conical conformal} this measure is non-atomic, ergodic and every other conformal measure of exponent $\alpha(f)$ is proportional to it.
\end{proof}

\bibliographystyle{plain}

\end{document}